\font\fifteenrm=cmr10 scaled\magstep2 
\font\fifteeni=cmmi10 scaled\magstep2
\font\fifteensy=cmsy10 scaled\magstep2
\font\fifteenbf=cmbx10 scaled\magstep2
\font\fifteentt=cmtt10 scaled\magstep2
\font\fifteenit=cmti10 scaled\magstep2
\font\fifteensl=cmsl10 scaled\magstep2
\font\fifteenam=msam10 scaled\magstep2
\font\fifteenbm=msbm10 scaled\magstep2
\font\fifteenex=cmex10 scaled\magstep2
\font\fifteensc=cmcsc10 scaled\magstep2 
\font\twelverm=cmr10 at 12pt
\font\twelvei=cmmi10 at 12pt
\font\twelvesy=cmsy10 at 12pt
\font\twelvebf=cmbx10 at 12pt
\font\twelvett=cmtt10 at 12pt
\font\twelveit=cmti10 at 12pt
\font\twelvesl=cmsl10 at 12pt
\font\twelveam=msam10 at 12pt
\font\twelvebm=msbm10 at 12pt
\font\twelveex=cmex10 at 12pt
\font\twelvesc=cmcsc10 at 12pt
\font\elevenrm=cmr10 scaled\magstephalf 
\font\eleveni=cmmi10 scaled\magstephalf
\font\elevensy=cmsy10 scaled\magstephalf
\font\elevenbf=cmbx10 scaled\magstephalf
\font\eleventt=cmtt10 scaled\magstephalf
\font\elevenit=cmti10 scaled\magstephalf
\font\elevensl=cmsl10 scaled\magstephalf
\font\elevenam=msam10 scaled\magstephalf
\font\elevenbm=msbm10 scaled\magstephalf
\font\elevenex=cmex10 scaled\magstephalf
\font\elevensc=cmcsc10 scaled\magstephalf
\font\tenrm=cmr10
\font\teni=cmmi10
\font\tensy=cmsy10
\font\tenbf=cmbx10
\font\tentt=cmtt10
\font\tenit=cmti10
\font\tensl=cmsl10
\font\tenam=msam10
\font\tenbm=msbm10
\font\tenex=cmex10
\font\tensc=cmcsc10
\font\ninerm=cmr9
\font\ninei=cmmi9
\font\ninesy=cmsy9
\font\ninebf=cmbx9
\font\ninett=cmtt9
\font\nineit=cmti9
\font\ninesl=cmsl9
\font\nineam=msam9
\font\ninebm=msbm9
\font\nineex=cmex9
\font\ninesc=cmcsc9
\font\eightrm=cmr8
\font\eighti=cmmi8
\font\eightsy=cmsy8
\font\eightbf=cmbx8
\font\eighttt=cmtt8
\font\eightit=cmti8
\font\eightsl=cmsl8
\font\eightam=msam8
\font\eightbm=msbm8
\font\eightex=cmex8
\font\eightsc=cmcsc8
\font\sevenrm=cmr7
\font\seveni=cmmi7
\font\sevensy=cmsy7
\font\sevenbf=cmbx7

\font\sevenam=msam7
\font\sevenbm=msbm7

\font\sixrm=cmr6
\font\sixi=cmmi6
\font\sixsy=cmsy6

\font\sixam=msam6
\font\sixbm=msbm6

\font\fiverm=cmr5
\font\fivei=cmmi5
\font\fivesy=cmsy5

\font\fiveam=msam5
\font\fivebm=msbm5

\font\fourrm=cmr5 at 4pt
\font\fouri=cmmi5 at 4pt
\font\foursy=cmsy5 at 4pt

\font\fouram=msam5 at 4pt
\font\fourbm=msbm5 at 4pt

\skewchar\twelvei='177 \skewchar\eleveni='177\skewchar\teni='177
\skewchar\ninei='177 \skewchar\eighti='177\skewchar\seveni='177 
\skewchar\sixi='177 \skewchar\fivei='177 \skewchar\fouri='177
\skewchar\twelvesy='60 \skewchar\elevensy='60 \skewchar\tensy='60
\skewchar\ninesy='60 \skewchar\eightsy='60 \skewchar\sevensy='60 
\skewchar\sixsy='60 \skewchar\fivesy='60 \skewchar\foursy='60
\newfam\itfam
\newfam\slfam
\newfam\bffam
\newfam\ttfam
\newfam\scfam
\newfam\amfam
\newfam\bmfam
\def\eightbig#1{{\hbox{$\left#1\vbox to 6.5pt{}\voidright $}}}
\def\eightBig#1{{\hbox{$\left#1\vbox to 7.5pt{}\voidright $}}}
\def\eightbigg#1{{\hbox{$\left#1\vbox to 10pt{}\voidright $}}}
\def\eightBigg#1{{\hbox{$\left#1\vbox to 13pt{}\voidright $}}}
\def\ninebig#1{{\hbox{$\left#1\vbox to 7.5pt{}\voidright $}}}
\def\nineBig#1{{\hbox{$\left#1\vbox to 8.5pt{}\voidright $}}}
\def\ninebigg#1{{\hbox{$\left#1\vbox to 11.5pt{}\voidright $}}}
\def\nineBigg#1{{\hbox{$\left#1\vbox to 14.5pt{}\voidright $}}}
\def\tenbig#1{{\hbox{$\left#1\vbox to 8.5pt{}\voidright $}}}
\def\tenBig#1{{\hbox{$\left#1\vbox to 9.5pt{}\voidright $}}}
\def\tenbigg#1{{\hbox{$\left#1\vbox to 12.5pt{}\voidright $}}}
\def\tenBigg#1{{\hbox{$\left#1\vbox to 16pt{}\voidright $}}}
\def\elevenbig#1{{\hbox{$\left#1\vbox to 9pt{}\voidright $}}}
\def\elevenBig#1{{\hbox{$\left#1\vbox to 10.5pt{}\voidright $}}}
\def\elevenbigg#1{{\hbox{$\left#1\vbox to 14pt{}\voidright $}}}
\def\elevenBigg#1{{\hbox{$\left#1\vbox to 17.5pt{}\voidright $}}}
\def\twelvebig#1{{\hbox{$\left#1\vbox to 10pt{}\voidright $}}}
\def\twelveBig#1{{\hbox{$\left#1\vbox to 11pt{}\voidright $}}}
\def\twelvebigg#1{{\hbox{$\left#1\vbox to 15pt{}\voidright $}}}
\def\twelveBigg#1{{\hbox{$\left#1\vbox to 19pt{}\voidright $}}}
\def\fifteenbig#1{{\hbox{$\left#1\vbox to 12pt{}\voidright $}}}
\def\fifteenBig#1{{\hbox{$\left#1\vbox to 13.5pt{}\voidright $}}}
\def\fifteenbigg#1{{\hbox{$\left#1\vbox to 18pt{}\voidright $}}}
\def\fifteenBigg#1{{\hbox{$\left#1\vbox to 23pt{}\voidright $}}}
\def\voidright{\right.\nulldelimiterspace=0pt \mathsurround=0pt }
\def\fifteenpoint{
  \textfont0=\fifteenrm \scriptfont0=\twelverm \scriptscriptfont0=\tenrm
  \def\rm{\fam0 \fifteenrm}%
  \textfont1=\fifteeni \scriptfont1=\twelvei \scriptscriptfont1=\teni
  \textfont2=\fifteensy \scriptfont2=\twelvesy \scriptscriptfont2=\tensy
  \textfont3=\fifteenex \scriptfont3=\fifteenex \scriptscriptfont3=\fifteenex
  \def\it{\fam\itfam\fifteenit}\textfont\itfam=\fifteenit
  \def\sl{\fam\slfam\fifteensl}\textfont\slfam=\fifteensl
  \def\bf{\fam\bffam\fifteenbf}\textfont\bffam=\fifteenbf 
    \scriptfont\bffam=\twelvebf\scriptscriptfont\bffam=\tenbf
  \def\tt{\fam\ttfam\fifteentt}\textfont\ttfam=\fifteentt
  \def\sc{\fam\scfam\fifteensc}\textfont\scfam=\fifteensc
  \def\am{\fam\amfam\fifteenam}\textfont\amfam=\fifteenam
    \scriptfont\amfam=\twelveam\scriptscriptfont\amfam=\tenam
  \def\bm{\fam\bmfam\fifteenbm}\textfont\bmfam=\fifteenbm
    \scriptfont\bmfam=\twelvebm\scriptscriptfont\bmfam=\tenbm
  \baselineskip=21pt \rm
  \let\big=\fifteenbig\let\Big=\fifteenBig\let\bigg=\fifteenbigg
  \let\Bigg=\fifteenBigg}
\def\twelvepoint{
  \textfont0=\twelverm \scriptfont0=\ninerm \scriptscriptfont0=\sevenrm
  \def\rm{\fam0 \twelverm}%
  \textfont1=\twelvei \scriptfont1=\ninei \scriptscriptfont1=\seveni
  \textfont2=\twelvesy \scriptfont2=\ninesy \scriptscriptfont2=\sevensy
  \textfont3=\twelveex \scriptfont3=\twelveex \scriptscriptfont3=\twelveex
  \def\it{\fam\itfam\twelveit}\textfont\itfam=\twelveit
  \def\sl{\fam\slfam\twelvesl}\textfont\slfam=\twelvesl
  \def\bf{\fam\bffam\twelvebf}\textfont\bffam=\twelvebf 
    \scriptfont\bffam=\ninebf\scriptscriptfont\bffam=\sevenbf
  \def\tt{\fam\ttfam\twelvett}\textfont\ttfam=\twelvett
  \def\sc{\fam\scfam\twelvesc}\textfont\scfam=\twelvesc
  \def\am{\fam\amfam\twelveam}\textfont\amfam=\twelveam
    \scriptfont\amfam=\nineam\scriptscriptfont\amfam=\sevenam
  \def\bm{\fam\bmfam\twelvebm}\textfont\bmfam=\twelvebm
    \scriptfont\bmfam=\ninebm\scriptscriptfont\bmfam=\sevenbm
  \baselineskip=17.8pt \rm 
  \def\looselineskip{\baselineskip=18.5pt plus 1.8pt}%
  \def\tightlineskip{\baselineskip=16.5pt}%
  \def\verytightlineskip{\baselineskip=15pt}%
  \let\big=\twelvebig\let\Big=\twelveBig\let\bigg=\twelvebigg
  \let\Bigg=\twelveBigg  }
\def\elevenpoint{
  \textfont0=\elevenrm \scriptfont0=\ninerm \scriptscriptfont0=\sixrm
  \def\rm{\fam0 \elevenrm}%
  \textfont1=\eleveni \scriptfont1=\ninei \scriptscriptfont1=\sixi
  \textfont2=\elevensy \scriptfont2=\ninesy \scriptfont2=\sixsy 
  \textfont3=\elevenex \scriptfont3=\elevenex \scriptfont3=\elevenex
  \def\it{\fam\itfam\elevenit}\textfont\itfam=\elevenit
  \def\sl{\fam\slfam\elevensl}\textfont\slfam=\elevensl
  \def\bf{\fam\bffam\elevenbf}\textfont\bffam=\elevenbf
  \def\tt{\fam\ttfam\eleventt}\textfont\ttfam=\eleventt
  \def\sc{\fam\scfam\elevensc}\textfont\scfam=\elevensc
  \def\am{\fam\amfam\elevenam}\textfont\amfam=\elevenam
    \scriptfont\amfam=\nineam\scriptscriptfont\amfam=\sixam
  \def\bm{\fam\bmfam\elevenbm}\textfont\bmfam=\elevenbm
    \scriptfont\bmfam=\ninebm\scriptscriptfont\bmfam=\sixbm
  \baselineskip=15.1pt \rm
  \def\looselineskip{\baselineskip=16pt plus 1.5pt}%
  \def\tightlineskip{\baselineskip=14pt}%
  \def\verytightlineskip{\baselineskip=13pt}%
  \let\big=\elevenbig\let\Big=\elevenBig\let\bigg=\elevenbigg
  \let\Bigg=\elevenBigg  }
\def\tenpoint{
  \textfont0=\tenrm \scriptfont0=\eightrm \scriptscriptfont0=\fiverm
  \def\rm{\fam0 \tenrm}%
  \textfont1=\teni \scriptfont1=\eighti \scriptscriptfont1=\fivei
  \textfont2=\tensy \scriptfont2=\eightsy \scriptfont2=\fivesy 
  \textfont3=\tenex \scriptfont3=\tenex \scriptfont3=\tenex
  \def\it{\fam\itfam\tenit}\textfont\itfam=\tenit
  \def\sl{\fam\slfam\tensl}\textfont\slfam=\tensl
  \def\bf{\fam\bffam\tenbf}\textfont\bffam=\tenbf
  \def\tt{\fam\ttfam\tentt}\textfont\ttfam=\tentt
  \def\sc{\fam\scfam\tensc}\textfont\scfam=\tensc
  \def\am{\fam\amfam\tenam}\textfont\amfam=\tenam
    \scriptfont\amfam=\eightam \scriptscriptfont\amfam=\fiveam
  \def\bm{\fam\bmfam\tenbm}\textfont\bmfam=\tenbm
    \scriptfont\bmfam=\eightbm \scriptscriptfont\bmfam=\fivebm
  \baselineskip=14pt\rm
  \def\looselineskip{\baselineskip=14.8pt plus1.5pt}
  \def\tightlineskip{\baselineskip=12.6pt}%
  \def\verytightlineskip{\baselineskip=13pt}%
  \let\big=\tenbig\let\Big=\tenBig\let\bigg=\tenbigg\let\Bigg=\tenBigg  }
\def\ninepoint{
  \textfont0=\ninerm \scriptfont0=\sevenrm \scriptscriptfont0=\fourrm
  \def\rm{\fam0 \ninerm}%
  \textfont1=\ninei \scriptfont1=\seveni \scriptscriptfont1=\fouri
  \textfont2=\ninesy \scriptfont2=\sevensy \scriptfont2=\foursy 
  \textfont3=\nineex \scriptfont3=\nineex \scriptfont3=\nineex
  \def\it{\fam\itfam\nineit}\textfont\itfam=\nineit
  \def\sl{\fam\slfam\ninesl}\textfont\slfam=\ninesl
  \def\bf{\fam\bffam\ninebf}\textfont\bffam=\ninebf
  \def\tt{\fam\ttfam\ninett}\textfont\ttfam=\ninett
  \def\sc{\fam\scfam\ninesc}\textfont\scfam=\ninesc
  \def\am{\fam\amfam\nineam}\textfont\amfam=\nineam
    \scriptfont\amfam=\nineam\scriptscriptfont\amfam=\fouram
  \def\bm{\fam\bmfam\ninebm}\textfont\bmfam=\ninebm
    \scriptfont\bmfam=\ninebm\scriptscriptfont\bmfam=\fourbm
  \baselineskip=12.6pt\rm
  \def\tightlineskip{\baselineskip=11.5pt}
  \let\big=\ninebig\let\Big=\nineBig\let\bigg=\ninebigg
  \let\Bigg=\nineBigg  }
\def\eightpoint{
  \textfont0=\eightrm \scriptfont0=\fiverm \scriptscriptfont0=\fourrm
  \def\rm{\fam0 \eightrm}%
  \textfont1=\eighti \scriptfont1=\fivei \scriptscriptfont1=\fouri
  \textfont2=\eightsy \scriptfont2=\fivesy \scriptfont2=\foursy 
  \textfont3=\eightex \scriptfont3=\eightex \scriptfont3=\eightex
  \def\it{\fam\itfam\eightit}\textfont\itfam=\eightit
  \def\sl{\fam\slfam\eightsl}\textfont\slfam=\eightsl
  \def\bf{\fam\bffam\eightbf}\textfont\bffam=\eightbf
  \def\tt{\fam\ttfam\eighttt}\textfont\ttfam=\eighttt
  \def\sc{\fam\scfam\eightsc}\textfont\scfam=\eightsc
  \def\am{\fam\amfam\eightam}\textfont\amfam=\eightam
    \scriptfont\amfam=\eightam\scriptscriptfont\amfam=\fouram
  \def\bm{\fam\bmfam\eightbm}\textfont\bmfam=\eightbm
    \scriptfont\bmfam=\eightbm\scriptscriptfont\bmfam=\fourbm
  \baselineskip=11.2pt \rm
  \def\tightlineskip{\baselineskip=10.4pt}
  \let\big=\eightbig\let\Big=\eightBig\let\bigg=\eightbigg
  \let\Bigg=\eightBigg  }

%
\twelvepoint
\nopagenumbers
\hsize=6in\vsize=8.8in

\parskip=1pt plus 1pt

\newif\ifSpecialhead\Specialheadfalse
\newbox\specialheadbox

\def\specialhead #1\par{\Specialheadtrue\setbox\specialheadbox=\hbox{#1}}
\headline={{\ifSpecialhead\box\specialheadbox\global\Specialheadfalse\else
     \ifnum\pageno<0{\hfill\quad{\twelvebf\folio}}%
     \else\ifnum\pageno<2\hfill
     \else\hfill\twelvepoint\sc\firstmark\quad{\twelvebf\folio}\fi\fi\fi}}

\def\title#1\par{\medskip{\def\cr{\hfil\par\hfil}\hfil\fifteenbf#1\hfil\par}\medskip}
\def\subtitle#1\par{\centerline{\fifteenrm #1}\medskip}
\def\author#1\par{\medskip{\def\cr{\hfil\par\hfil\twelvesc}\fifteensc\hfil#1\hfil\par}}
\def\authors#1\par{\medskip\fifteensc\center#1\par}
\def\center#1\par{{\def\cr{\hfil\par\hfil}\hfil#1\hfil\par}}
\def\abstract.#1\par{\message{Abstract.}%
                    \medskip{\narrower\narrower\tenpoint\tightlineskip
                        \noindent{\bf Abstract.}#1\par}\medskip\noindent}
\def\tinyabstract.#1\par{\message{Abstract.}%
                    \medskip{\narrower\narrower\eightpoint\tightlineskip
                        \noindent{\bf Abstract.}#1\par}\medskip\noindent}
\def\bigabstract.#1\par{\message{Abstract.}%
                         \medskip{\narrower\narrower\tightlineskip
                         \noindent{\bf Abstract. }#1\par}\medskip\noindent}
\def\acknowledgement#1\par{\footnote{}{#1}}
\def\sectionskip{\Goodbreak\vskip 25pt plus 15pt minus 5pt}
\def\secnumber{\ifquiet
               \else\ifNoSections
                    \else\sectionsymbol\the\secno\quad\fi\fi}
\def\section#1\par{ \NoSectionsfalse\par\sectionskip\proofdepth=0\claimno=0
 \ifquiet\else\advance\secno by1\fi\toks0={#1}
 \immediate\write16{\ifquiet\else Section \the\secno\space\fi
                    \the\toks0}%
 \mark{\secnumber #1}%
 {\fifteenpoint\bf\noindent\secnumber #1}\nobreak\bigskip\quietoff
 \nobreak\noindent}
\def\quiet{\quiettrue}

\def\quietoff{\ifQUIET\else\quietfalse\fi}
\newif\ifquiet
\newif\ifQUIET
\newif\ifNoSections
\newcount\claimtype
\newcount\secno
\newcount\claimno
\newcount\subclaimno
\newcount\subsubclaimno
\newcount\subsubsubclaimno
\newcount\proofdepth
\def\subclaimnumber{\ifquiet\else\ifcase\subclaimno\or A\or B\or C\or D\or E\or
     F\or G\or H\or I\or J\or K\or L\or M\or N\or O\or P\fi\fi}
\def\subsubclaimnumber{\ifquiet\else\ifcase\subsubclaimno\or i\or ii\or iii\or
   iv\or v\or vi\or vii\or viii\or ix\or x\or xi\or xii\or xiii\or xiv\fi\fi}
\def\subsubsubclaimnumber{\ifquiet\else\ifcase\subsubsubclaimno\or a\or b\or
   c\or d\or e\or f\or g\or \or h\or i\or j\or k\or l\or m\or n\or o\fi\fi}
\def\claimtag{\ifquiet\else
  \ifNoSections
    \ifcase\proofdepth\the\claimno%
    \or\the\claimno.\subclaimnumber
    \or\the\claimno.\subclaimnumber.\subsubclaimnumber
    \or\the\claimno.\subclaimnumber.\subsubclaimnumber
                                                .\subsubsubclaimnumber\fi
  \else
    \ifcase\proofdepth\the\secno.\the\claimno
    \or\the\secno.\the\claimno.\subclaimnumber
    \or\the\secno.\the\claimno.\subclaimnumber.\subsubclaimnumber
    \or\the\secno.\the\claimno.\subclaimnumber.\subsubclaimnumber
                                                .\subsubsubclaimnumber\fi\fi\fi}
\secno=0\claimno=0\proofdepth=0\subclaimno=0\subsubclaimno=0\subsubsubclaimno=0
\NoSectionstrue
\newbox\qedbox
\def\claimname{\ifcase\claimtype Theorem\or Lemma\or Claim\or Corollary\or
               Question\or Definition\or Remark\or Conjecture\fi}
\def\preclaimskip{\removelastskip
    \ifcase\claimtype\goodbreak\vskip 8pt plus 4pt minus 2pt
                  \or\goodbreak\vskip 6pt plus 4pt minus 1pt
                  \or\goodbreak\vskip 5pt plus 4pt minus 1pt
                  \or\goodbreak\vskip 8pt plus 4pt minus 2pt
                  \or\vskip 7pt plus 4pt minus 2pt
                  \or\vskip 7pt plus 4pt minus 2pt
                  \or\vskip 7pt plus 4pt minus 2pt
                  \or\goodbreak\vskip 8pt plus 4pt minus 2pt\fi}
\def\postclaimskip{\ifcase\claimtype         \vskip 4pt plus 2pt minus 2pt
                                          \or\vskip 3pt plus 2pt minus 2pt
                                          \or\vskip 2pt plus 2pt minus 1pt
                                          \or\vskip 4pt plus 2pt minus 2pt
                                          \or\vskip 1pt plus 2pt
                                          \or\vskip 4pt plus 4pt
                                          \or\vskip 3pt plus 2pt
                                          \or\vskip 4pt plus 2pt minus 2pt\fi}
\def\claimfont{\ifcase\claimtype
                  \sl\or\sl\or\sl\or\sl\or\sl\or\rm\or\rm\or\sl\fi}
\def\advancetag{\ifcase\proofdepth\advance\claimno by1
                               \or\advance\subclaimno by1
                               \or\advance\subsubclaimno by1
                               \or\advance\subsubsubclaimno by1\fi}
\def\sayclaim#1.#2 #3\par{\ifquiet\else\advancetag\fi
    \preclaimskip\setbox1=\hbox{#1}\setbox2=\hbox{#2}%
    \toks0={#1 }
    \immediate\write16{\ifdim\wd1>0pt\the\toks0
                       \else\claimname\space\fi \claimtag.}%
    \vbox{\noindent
    {\bf\ifdim\wd1=0pt \claimname\else #1\fi\ifquiet.\else\ \claimtag{\ifNoSections.\fi}\fi}%
    \enspace{\ifdim\wd2>0pt\sc #2\enspace\fi}%
    {\claimfont #3\par}}\postclaimskip\quietoff}
\def\theorem{\claimtype=0\sayclaim}
\def\lemma{\claimtype=1\sayclaim}

\def\corollary{\claimtype=3\sayclaim}
\def\question{\claimtype=4\sayclaim}
\def\definition{\claimtype=5\sayclaim}
\def\remark{\claimtype=6\sayclaim}
\def\conjecture{\claimtype=7\sayclaim}
\def\point#1. #2\par{\item{\rm #1.}#2\par}
\def\points#1\cr\par{\medskip\vbox{\let\cr=\point\point#1\par}\par}
\def\df{\it}
\def\prooffont{}
\def\proofsize{}
\def\proofindent{}
\def\proofskip{\badbreak\ifcase\claimtype    \vskip 3pt plus 2pt minus 2pt
                                          \or\vskip 2pt plus 2pt minus 2pt
                                          \or\vskip 1pt plus 2pt minus 1pt
                                          \or\vskip 3pt plus 2pt minus 2pt
                                          \or\vskip 1pt plus 2pt
                                          \or\vskip 2pt plus 4pt
                                          \or\vskip 1pt plus 2pt
                                          \or\vskip 3pt plus 2pt minus 2pt\fi}

\def\Goodbreak{\vskip0pt plus.5in\penalty-1000\vskip0pt plus-.5in}
\def\goodbreak{\penalty-500}
\def\badbreak{\penalty500}
\def\Badbreak{\penalty1000}
\def\proof{\message{proof}\removelastskip\Badbreak\proofskip\begingroup
  \advance\proofdepth by1
  \setbox\qedbox=\hbox{\halmos\raise2pt\hbox{\fiverm\claimname}}%
  \prooffont\proofsize\proofindent\noindent{\bf Proof: }}
\def\proofof#1:{\message{proof}\removelastskip\Badbreak\proofskip\begingroup
  \advance\proofdepth by1
  \setbox\qedbox=\hbox{\halmos\raise2pt\hbox{\fiverm#1}}%
  \prooffont\proofsize\proofindent\noindent{\bf Proof of #1: }}
\def\cite[#1]{[{\tenrm{#1}}]\message{[#1]}}
\edef\ref#1{\expandafter\global\expandafter\edef#1{\noexpand\claimtag}}
\newwrite\notes
\openout\notes=\jobname.notes
\long\def\unexpandedwrite#1#2{\def\finwrite{\write#1}%
   {\aftergroup\finwrite\aftergroup{\sanitize#2\endsanity}}}
\def\sanitize{\futurelet\next\sanswitch}
\let\stoken=\space
\def\sanswitch{\ifx\next\endsanity
  \else\ifcat\noexpand\next\stoken\aftergroup\space\let\next=\eat
   \else\ifcat\noexpand\next\bgroup\aftergroup{\let\next=\eat
    \else\ifcat\noexpand\next\egroup\aftergroup}\let\next=\eat
     \else\let\next=\copytoken\fi\fi\fi\fi \next}
\def\eat{\afterassignment\sanitize \let\next= }
\long\def\copytoken#1{\ifcat\noexpand#1\relax\aftergroup\noexpand
  \else\ifcat\noexpand#1\noexpand~\aftergroup\noexpand\fi\fi
  \aftergroup#1\sanitize}
\def\endsanity\endsanity{}

\def\note#1#2{\hbox to2in{\strut#1\quad\dotfill\quad#2}}
\def\boxit#1{\setbox4=\hbox{\kern1pt#1\kern1pt}
  \hbox{\vrule\vbox{\hrule\kern1pt\box4\kern1pt\hrule}\vrule}}
\def\halmos{\hbox{\am\char'3}}
\def\qed#1\par{\message{.                                }\setbox1=\hbox{#1}%
  \ifdim\wd1>0pt\setbox\qedbox=\hbox{\halmos\raise2pt\hbox{\fiverm #1}}\fi
  \kern5pt\lower 2pt\hbox{\box\qedbox}\proofskip\goodbreak\endgroup}

\def\sectionsymbol{\S}
\def\k{\kappa}
\def\g{\gamma}
\def\a{\alpha}
\def\b{\beta}
\def\d{\delta}
\def\s{\sigma}
\def\t{\tau}
\def\l{\lambda}
\def\z{\zeta}
\def\I1{\mathop{\hbox{\sc i}_1}}
\def\w{\omega}
\def\P{{\mathchoice{\hbox{\bm P}}{\hbox{\bm P}}
         {\hbox{\tenbm P}}{\hbox{\sevenbm P}}}}
\def\Q{{\mathchoice{\hbox{\bm Q}}{\hbox{\bm Q}}
         {\hbox{\tenbm Q}}{\hbox{\sevenbm Q}}}}
\def\R{{\mathchoice{\hbox{\bm R}}{\hbox{\bm R}}
         {\hbox{\tenbm R}}{\hbox{\sevenbm R}}}}

\def\card#1{\left|#1\right|}

\def\dom{\mathop{\rm dom}\nolimits}

\def\RO{\mathop{\hbox{\sc ro}}\nolimits}
\def\TC{\mathop{\hbox{\sc tc}}\nolimits}

\def\id{\mathop{\hbox{\tenrm id}}}

\def\elesub{\prec}

\def\unifto{\buildrel\lower 7pt\hbox{$\to$}\over\to}

\def\iso{\cong}

\def\cof{\mathop{\rm cof}\nolimits}
\def\cp{\mathop{\rm cp}\nolimits}
\def\ran{\mathop{\rm ran}\nolimits}
\def\from{\mathbin{\vbox{\baselineskip=3pt\lineskiplimit=0pt
                         \hbox{.}\hbox{.}\hbox{.}}}}
\def\ORD{\hbox{\sc ord}}

\def\ZFC{\hbox{\sc zfc}}
\def\GCH{\hbox{\sc gch}}
\def\SCH{\hbox{\sc sch}}

\def\plus{^{\scriptscriptstyle +}}
\def\plusplus{^{\scriptscriptstyle ++}}
\def\plusplusplus{^{\scriptscriptstyle +++}}
\def\plusplusplusplus{^{\scriptstyle ++++}}
\def\in{\mathrel{\mathchoice{\raise
1pt\hbox{$\scriptstyle\cal\char'62$}}
         {\raise 1pt\hbox{$\scriptstyle\cal\char'62$}}
         {\raise .5pt\hbox{$\scriptscriptstyle\cal\char'62$}}
         {\hbox{$\scriptscriptstyle\cal\char'62$}}}\penalty700{}}
\def\ni{\mathrel{\mathchoice{\raise 1pt\hbox{$\scriptstyle\cal\char'63$}}
                   {\raise 1pt\hbox{$\scriptstyle\cal\char'63$}}
                   {\raise .5pt\hbox{$\scriptscriptstyle\cal\char'63$}}
                   {\hbox{$\scriptscriptstyle\cal\char'63$}}}\penalty700}
\def\of{\mathrel{\mathchoice{\raise 1pt\hbox{$\scriptstyle\subseteq$}}
                   {\raise 1pt\hbox{$\scriptstyle\subseteq$}}
                   {\raise .5pt\hbox{$\scriptscriptstyle\subseteq$}}
                   {\hbox{$\scriptscriptstyle\subseteq$}}}}
\def\fo{\mathrel{\mathchoice{\raise 1pt\hbox{$\scriptstyle\supseteq$}}
                   {\raise 1pt\hbox{$\scriptstyle\supseteq$}}
                   {\raise .5pt\hbox{$\scriptscriptstyle\supseteq$}}
                   {\hbox{$\scriptscriptstyle\supseteq$}}}}
\def\notin{\mathrel{\mathchoice
  {\raise 1pt\hbox{\rlap{$\scriptstyle\;|$}$\scriptstyle\cal\char'62$}}
  {\raise 1pt\hbox{\rlap{$\scriptstyle\kern2pt
          |$}$\scriptstyle\cal\char'62$}}
  {\raise .5pt\hbox{\rlap{$\scriptscriptstyle\, |$}$\scriptscriptstyle
      \cal\char'62$}}
  {\hbox{\rlap{$\scriptscriptstyle\, |$}$\scriptscriptstyle
     \cal\char'62$}}}%
  \penalty700}

\def\and{\mathrel{\kern1pt\&\kern1pt}}
\def\iff{\mathrel{\leftrightarrow}}

\def\implies{\rightarrow}

\def\union{\cup}

\def\compose{\circ}
\def\intersect{\cap}

\def\ot{\mathop{\rm ot}\nolimits}
\def\setminus{\mathbin{\hbox{\bm\char'162}}}
\def\minus{\setminus}
\def\add{\mathop{\rm Add}\nolimits}

\def\cross{\times}

\def\muchgt{\gg}
\def\lt{\mathrel{\mathchoice{\scriptstyle<}{\scriptstyle<}
   {\scriptscriptstyle<}{\scriptscriptstyle<}}}
\def\lte{\mathrel{\scriptstyle\leq}}
\def\tlt{\triangleleft}

\def\[#1]{\left[\vphantom{\bigm|}#1\right]}
\def\<#1>{\langle\,#1\,\rangle}

\def\sat{\models}
\def\satisfies{\models}

\def\image{\mathbin{\hbox{\tt\char'42}}}
\def\restrict{\mathbin{\mathchoice{\hbox{\am\char'26}}{\hbox{\am\char'26}}{\hbox{\eightam\char'26}}{\hbox{\sixam\char'26}}}}
\def\force{\mathbin{\hbox{\am\char'15}}}

\def\emptyset{\mathord{\hbox{\bm\char'77}}}
\def\beth{\mathord{\hbox{\bm\char'151}}}
\def\boolval#1{\mathopen{\lbrack\!\lbrack}\,#1\,\mathclose{\rbrack\!\rbrack}}

\def\st{\mid}
\def\seq<#1>{{\def\st{\mid\penalty650}\left<\,#1\,\right>}}

\def\set#1{\{\,{#1}\,\}}

\def\th{{\hbox{\fiverm th}}}

\def\forces{\force}
\def\lttheta{{\raise 1pt\hbox{$\scriptstyle<$}\theta}}

\def\I1{\mathop{\hbox{\sc i}_1}}
\def\ltk{{{\scriptstyle<}\k}}

\def\ltg{{{\scriptstyle<}\g}}

\def\ltb{{{\scriptstyle<}\b}}
\def\lteb{{{\scriptstyle\leq}\b}}
\def\lted{{{\scriptstyle\leq}\d}}

\def\Qdot{\dot\Q}

\def\Qdot{\dot\Q}

\def\Pforces{\force_{\P}}

\def\jmu{j_\mu}

\def\Ptail{\P_{\fiverm \!tail}}

\def\ltg{{\scriptscriptstyle<}\g}

\def\Gtail{G_{\fiverm tail}}

\def\ltek{{\scriptstyle\leq}\k}

\def\Qtilde{{\widetilde \Q}}

\font\arrow=line10 scaled \magstep1
\def\makeline#1.{\hbox{\arrow\char#1}}
\def\makearrow#1.#2.{\hbox{\arrow\char#1\llap{\char#2}}}
\def\definelinesandarrows#1.#2.#3.#4.#5.{
   \expandafter\edef\csname#4line\endcsname{\makeline#1.}
   \expandafter\edef\csname#4arrow\endcsname{\makearrow#1.#2.}
   \expandafter\edef\csname#5line\endcsname{\makeline#1.}
   \expandafter\edef\csname#5arrow\endcsname{\makearrow#1.#3.}}
\definelinesandarrows 0.18.9.ne.sw.
\definelinesandarrows 1.21.11.nnne.sssw.
\definelinesandarrows 2.14.13.nnnne.ssssw.
\definelinesandarrows 3.23.15.nnnnne.sssssw.
\definelinesandarrows 4.23.15.nnnnnne.ssssssw.
\definelinesandarrows 10.30.29.nne.ssw.
\definelinesandarrows 16.49.41.neeeeee.swwwwww.
\definelinesandarrows 17.51.43.neeee.swwww.
\definelinesandarrows 19.55.47.nehuh.swhuh.
\definelinesandarrows 24.58.41.neeeeeee.swwwwwww.
\definelinesandarrows 26.62.9.neee.swww.
\definelinesandarrows 33.49.25.neeeee.swwwww.
\definelinesandarrows 35.62.61.nee.sww.
\definelinesandarrows 64.82.73.se.nw.
\definelinesandarrows 65.85.75.ssse.nnnw.
\definelinesandarrows 66.78.77.sssse.nnnnw.
\definelinesandarrows 67.87.79.ssssse.nnnnnw.
\definelinesandarrows 68.87.79.sssssse.nnnnnnw.
\definelinesandarrows 74.94.93.sse.nnw.
\definelinesandarrows 80.113.105.seeeeee.nwwwwww.
\definelinesandarrows 81.115.107.seeee.nwwww.
\definelinesandarrows 99.126.125.see.nww.
\def\sejoin#1#2{\setbox1=\hbox{#1}\setbox2=\hbox{#2}%
  \hbox{\vbox{\hbox{\copy1\kern\wd2}\nointerlineskip
              \hbox{\kern\wd1\box2}}}}
\def\nejoin#1#2{\setbox1=\hbox{#1}\setbox2=\hbox{#2}%
  \hbox{\vbox{\hbox{\kern\wd1\copy2}\nointerlineskip\hbox{\copy1\kern\wd2}}}}
\newdimen\hnudge
\newdimen\vnudge
\newdimen\hnudgedefault
\newdimen\vnudgedefault

\def\SEdefaultnudge{\hnudge=-16pt\vnudge=20pt}
\def\Edefaultnudge{\hnudge=-25pt\vnudge=6pt}
\def\Sdefaultnudge{\hnudge=-8pt\vnudge=20pt}
\def\longEdefaultnudge{\hnudge=-5pt\vnudge=6pt}
\def\nudgeright#1pt{\advance\hnudge by#1pt}
\def\nudgeleft#1pt{\advance\hnudge by-#1pt}
\def\nudgeup#1pt{\advance\vnudge by#1pt}
\def\nudgedown#1pt{\advance\vnudge by-#1pt}
\def\label#1{\smash{\llap{\kern\hnudge
                   \raise\vnudge\rlap{$\scriptstyle#1$}\hfill}}}

\def\SEarrow{\SEdefaultnudge
             \sejoin\seeline{\sejoin\seeline{\sejoin\seeline\seearrow}}}

\def\Sarrow{\Sdefaultnudge\setbox1=\hbox{\SEarrow}
           \hbox{\hskip 10pt\vrule height\ht1\hbox{\arrow\char'77}}}
\def\Earrow{\Edefaultnudge\setbox1=\hbox{\SEarrow}
 \hbox{\raise 2pt\hbox{\vrule height-.4pt depth.8ptwidth\wd1\kern2pt
       \llap{\arrow\char'55}}}}
\def\longEarrow{\longEdefaultnudge\setbox1=\hbox{\SEarrow}
      \rlap{\hskip-1.25\wd1\raise 2pt
            \hbox{\vrule height-.4pt depth.8ptwidth2.5\wd1\kern2pt
            \llap{\arrow\char'55}}}}
\def\trianglediagram#1#2#3#4#5#6{%
    {\def\normalbaselines{\baselineskip0pt\lineskip8pt\lineskiplimit0pt}%
       \matrix{#1& &\cr
               \Sarrow\label{#2}&\SEarrow\label{#3}&\cr
               #4&\Earrow\label{#5}&#6\cr}}}

\def\Ptail{\P_{\!\!\hbox{\fiverm tail}}}
\def\Gtail{G_{\!\hbox{\fiverm tail}}}
\def\gtail{g_{\!\hbox{\fiverm tail}}}
\def\Qtilde{\tilde\Q}
\def\ltg{{\lt}\g}

\def\ltek{{\leq}\k}

\title Destruction or Preservation As You Like It

\author Joel David Hamkins

\abstract. The Gap Forcing Theorem, a key contribution of this paper,
implies essentially that after any reverse Easton iteration of closed
forcing, such as the Laver preparation, every supercompactness measure
on a supercompact cardinal extends a measure from the ground model.
Thus, such forcing can create no new supercompact cardinals, and, if the
GCH holds, neither can it increase the degree of supercompactness of
any cardinal; in particular, it can create no new measurable cardinals.
In a crescendo of what I call exact preservation
theorems, I use this new technology to perform a kind of partial Laver
preparation, and thereby finely control the class of
posets which preserve a supercompact cardinal. Eventually, I prove the
`As You Like It' Theorem, which asserts that the class of
$\ltk$-directed closed posets which preserve a supercompact cardinal
$\k$ can be made by forcing to conform with any pre-selected local
definition which respects the equivalence of forcing. Along the way I
separate completely the levels of the superdestructibility hierarchy,
and, in an epilogue, prove that the notions of fragility and
superdestructibility are orthogonal---all four combinations are
possible.

There is a vast unknown continent, which I aim to explore, between two
extreme theorems:  Laver \cite[Lav78], on the one hand, proved, in what is perhaps
my favorite argument in large cardinal set theory, that a supercompact
cardinal can be made {\it indestructible}, so that its supercompactness
is preserved by every $\ltk$-directed closed forcing notion;
in my Superdestruction Theorem \cite[Ham97b], on the other hand, and in joint work
with Shelah \cite[HS97], I proved that a supercompact cardinal can be made {\it
superdestructible}, so that its supercompactness  is destroyed by every
$\ltk$-closed forcing notion. Are there any
theorems in the uncharted wilderness between these two extremes?
Indeed there are, and in this paper I will prove that there are. Here,
between indestructibility and superdestructibility, I will lean
alternately on the methods of \cite[Lav78] when I want a poset to
preserve supercompactness and then on the Gap Forcing Theorem,
introduced in this paper, when I want it to destroy supercompactness.
In a crescendo of Exact Preservation theorems, my results will
culminate in the `As You Like It' Theorem, which asserts that the class
of $\ltk$-directed closed posets which preserve the supercompactness of
$\k$ can be made to conform with any pre-selected local definition
which respects the equivalence of forcing. This theorem and the others
I prove like it fill the region between indestructibility and
superdestructibility. I hope that my techniques will allow you to prove
that the class of posets which preserve the supercompactness of $\k$
can be made to be the class you have always wished that it was,
whatever that may be.

I confess that this paper is inside-out. After proving some preliminary
facts, I march through a sequence of exact preservation theorems, the
later ones stronger versions of earlier ones, finally advancing to the
As You Like It Theorem.  In an outside-out fashion, I could simply have
proved the As You Like It Theorem first, and then deduced the earlier
theorems as corollaries. But I chose this backwards manner of
presentation because I view the main contribution of this paper as the
{\it method of proof} of these theorems. So, in order to highlight the
power of this method, I start slowly and then build up to the stronger
theorems. I hope that the logical overlap that this order of
presentation involves will be forgiven.

Let me explain the paper's overall structure. I begin in section one
with my main new tool, the Gap Forcing Theorem. In sections two and
three I introduce the useful concepts of a partial Laver preparation
and a high jump function, respectively, before giving my first
applications in section four: the Exact Preservation Theorems. In
section five, in order to improve these theorems, I develop the theory
of representability, and, in section six, apply it to separate the
levels of the superdestructibility hierarchy. In section seven, I
present much more powerful Exact Preservation Theorems, and culminate
in the `As You Like It' Theorem, the title theorem of this paper.
Finally, in an epilogue, I separate the notions of superdestructibility
and fragility.

I try whenever possible to use standard notation, but assume a
familiarity with reverse Easton forcing iterations, such as the Laver
preparation, and the lifting arguments they typically involve.
Following Adrian Mathias, I use the notation $h\from X\to Y$ to mean
that $h$ is a partial function from $X$ to $Y$.

\section Gap Forcing

The key new technology in this paper is the Gap Forcing Theorem and its
corollaries, which give explicit information about the
nature of supercompactness embeddings which are added by forcing. Since they
severely limit the kinds of measures which can exist after gap forcing,
I will use them as a fundamental tool when proving the Exact
Preservation Theorems, where I must make certain forcing notions
destroy supercompactness.

The Gap Forcing Theorem will depend crucially on an improved version of
The Key Lemma, below, which was proved initially in \cite[Ham97b] but was
also instrumental in the main results of \cite[HamShl].  As in \cite[Ham97b],
define that a sequence or a set of ordinals is {\df fresh} over $V$
when it is not in $V$, but every proper initial segment of it is in
$V$.

\lemma Key Lemma. If $\card{\P}=\b$, $\cof(\l)>\b$, and $\Pforces\Q$ is
$\lteb$-closed, then $\P*\Q$ adds no fresh subsets of $\l$.
\ref\KeyLemma

In this paper I will use an improved version of the Key Lemma. Before doing
so, let me define that a poset is {\it $\lteb$-strategically closed} when
the second player has a winning strategy in the game of length $\b$ in
which the players create a descending sequence $\<p_\a\st \a<\b>$ from
the poset, the second player playing at every limit stage. The first
player to violate the rule that the conditions descend loses, and
otherwise the second player wins.

\lemma Improved Key Lemma. If $\card{\P}=\b$, $\cof(\l)>\b$, and
$\Pforces\Q$ is $\lteb$-strategically closed, then $\P*\Q$ adds no fresh
subsets of $\l$, and no fresh $\l$-sequences.
\ref\ImprovedKeyLemma

\proof There are two improvements. First, the original proof of the Key
Lemma in \cite[Ham97b] works in the case, here, of strategically closed
$\Q$, just as well: one simply needs to take care, when defining $q_t$
in that argument, to also obey the strategy, so that the limit stages
will go through. Second, the improvement to $\l$-sequences is obtained
by coding elements of $\d^\l$ with binary sequences of length
$\d\cdot\l$, which has the same cofinality as $\l$.\qed

\definition Main Definition. A forcing notion $\P$ admits a {\df gap}
at $\d$ when the poset $\P$ factors as $\P_1*\P_2$ where
$\card{\P_1}<\d$ and $\P_2$ is $\lted$-strategically closed in
$V^{\P_1}$. The Laver preparation, and indeed most every reverse Easton
iteration of closed forcing, admits a gap between any two stages of the
forcing.

\definition. A set $C$ is {\df unbounded} in $P_\k\g$ when for every
$\s\in P_\k\g$ there is $\t\in C$ with $\s\of \t$. A set $D\of C$ is
{\df $\d$-directed} if for any set $B\of D$ of size less than $\d$
there is $\s\in D$ such that $\union B\of\s$.  The set $C$ is {\df
$\d$-closed} when every $\d$-directed $D\of C$ with $\card{D}<\k$ has
$\union D\in C$.  Define that $C$ is {\df $\d$-club} when it is
unbounded and $\d$-closed (Caution: this usage differs from that in,
say, \cite[Kan94]).  A {\df supercompactness} measure is simply a
normal fine measure on $P_\k\l$ for some $\l$.

I view the next theorem, evolved from lesser forms into the highly
useful current animal, as a full-grown version of the Key Lemma. It and
its corollaries will provide the essential new tools with which I will
deduce that certain kinds of forcing notions destroy supercompactness.

You may underappreciate this theorem if you are not familiar with the
bizarre sorts of embeddings which can live in a forcing extension. It
is easy, for example, to construct a forcing extension $V[G]$ with an
embedding $j:V[G]\to M[j(G)]$ such that $M\not\of V$. Indeed, $M$ may
even have different subsets of $\k=\cp(j)$ than $V$.  Adding a Cohen
subset to a supercompact cardinal $\k$ which is indestructible in $V$,
for example, of necessity produces such embeddings, since the new
subset $G\of\k$ must be in $M[j(G)]$ but cannot have been added by
$j(G)$, and so it is in $M$ but not in $V$.  Of course, however, a
standard technique to show supercompactness is preserved to a forcing
extension is to lift an embedding $j:V\to M$ from the ground model $V$
to the extension $j:V[G]\to M[j(G)]$, and these embeddings have the
properties listed in the conclusion of the Gap Forcing Theorem.  But
that is exactly the force and utility of the theorem---the amazing fact
the Gap Forcing Theorem and its corollaries identify is that after
forcing which admits a gap below $\k$, {\it every} supercompactness
embedding resembles the lift of a supercompactness embedding in $V$.
But this is only resemblance:  we cannot prove that every
$\l$-supercompactness measure in a forcing extension $V[G]$ which
admits a gap below $\k$ is the lift of a $\l$-supercompactness measure
in $V$. This is because, as I proved in my dissertation \cite[Ham94b], a
strong embedding in $V$, which is not an ultrapower embedding at all,
can be lifted via gap forcing to an embedding which in the extension is
the ultrapower by a normal measure on $\k$.  So we should perhaps be
content to know that the embeddings in $V[G]$ resemble lifts as much as
the claims made in the Gap Forcing Theorem require them to. Of course,
I am speaking here of embeddings which are internal to $V[G]$ in the
sense that they are definable there.

\theorem Gap Forcing Theorem. Suppose that $V[G]$ is a forcing
extension of $V$ which admits a gap below $\k$, and that $j:V[G]\to
M[j(G)]$ is an embedding defined in $V[G]$, with critical point $\k$,
which is closed under $\l$-sequences, i.e. $M[j(G)]^\l\of M[j(G)]$ in
$V[G]$, for some $\l$ with $\k\leq\l<j(\k)$. Then:
\points 1. $M\of V$,\cr
        2. $(M^\l)^M=(M^\l)^V$,\cr
        3. $j\image\l\in M$, and\cr
        4. $j\restrict\l^{+^{(n)V}}\in V$ for every $n\in\w$. (This can
            be pushed higher depending on the distributivity of the forcing.)\cr

\ref\GapForcingTheorem

\proof Suppose that $G=g*H\of\P*\Q$ exhibits the gap at $\d<\k$, so
that $\card{\P}<\d$ and $\Q$ is $\lted$-strategically closed in $V[g]$.
Without loss of generality I may assume that $\P\in V_\k$, so that
$j(\P)=\P$.

Before proving the claims, let me first establish the preliminary claim that
$M^\l\intersect V=V^\l\intersect M$. It suffices to show by
induction on $\g\le\l$ that $M^\g\intersect V=V^\g\intersect M$.
Suppose this is true up to $\g$, and consider now the two directions to
prove at $\g$. Suppose, for $(\of)$, that $s\in M^\g\intersect V$.
Thus, $s:\g\to M$ and $s\in V$; I aim to show $s\in M$. By the
induction hypothesis every initial segment of $s$ is in $M$. Certainly
$s\in M[j(G)]$, by the closure assumption on $j$.
If $\cof(\g)>\d$, then $s$ must in fact be in $M$, for otherwise it
would be an $M$-fresh sequence added by the gap forcing $j(\P*\Q)$,
contrary to the Key Lemma \ImprovedKeyLemma. So assume that $\cof(\g)\le\d$. Thus, by the
closure of $j(\Q)$, it follows that $s\in M[g]$, and so $s=\dot s_g$
for some name $\dot s\in M$. View $\dot s$ as a function from $\g$ to
the set of antichains in $\P$ labelled with the possible elements of
$\ran(s)$, which are all in both $M$ and $V$. By the closure of $\Q$,
it follows that $\dot s\in V[g]$, and so $\dot s=\ddot s_g$ for some
name-of-a-name $\ddot s\in V$. In $V[g]$ we have $(\ddot s_g)_g=s$, and
so this is forced by some condition $p\in g$.  The condition $p$ must
decide enough information about the name-of-a-name $\ddot s$ to
determine that it is a name which $p$ decides agrees with the sequence
$s$. (This is an unusual use of a name-of-a-name in that unlike iterated
forcing, here I am using the {\it same} generic twice to interpret
the name in $(\ddot s_g)_g$.)
Whatever information $p$ decides about $\ddot s$ must agree
with $\dot s$, since $\dot s=\ddot s_g$, so it follows that from the
condition $p$ and $\dot s$ in $M$ I can construct $s$. So $s\in M$.

For the converse direction, suppose $s:\g\to V$ and $s\in M$. I aim to
show that $s\in V$. By the induction hypothesis, every initial segment
of $s$ is in $V$, and so if $\cof(\g)>\d$, the Key Lemma \ImprovedKeyLemma\ yields $s\in
V$. So assume again that $\cof(\g)\le\d$. Thus $s\in V[g]$ by the
closure of $\Q$, and consequently $s=\dot s_g$ for some name $\dot s\in
V$.  By coding the name $\dot s$ as a $\g$-sequence, it follows by the
previous direction that $\dot s\in M$. Thus, both $s$ and $\dot s$ are
in $M$, and in $M[g]$ we have $s=\dot s_g$. Thus, there is a condition
$p\in g$ forcing this. Using $p$ and $\dot s$ in $V$, I can construct
$s$.  So $s\in V$, as desired. This proves the preliminary claim.

Now let me prove 1. It suffices to show that $P(\b)^M\of V$ for every
ordinal $\b$. Suppose $X\in P(\b)^M$ and, by induction, every initial
segment of $X$ is in $V$. I aim to show $X\in V$. If $\cof(\b)>\d$,
then $X$ is in $V$ by the Key Lemma \KeyLemma. So assume, alternatively, that
$\cof(\b)\le\d$, and consequently that $\b=\sup_{\a<\d}\b_\a$ in $M$, with
$\b_\a<\b$. In this case $\<X\intersect\b_\a\st\a<\d>$ is in
$V^\d\intersect M$, and consequently it is in $V$ by the preliminary
claim. By taking a
union it follows that $X\in V$ as well. This establishes 1. And 2
follows immediately from 1 and the preliminary claim.

To prove 3, it suffices to show by induction on $\g$
that $j\image\g\in M$ for all $\g\le\l$. Assume by induction that this
claim holds up to $\g$, and consider $j\image\g$. I may assume
$\g>\k$.  Note that $j\image\l\in M[j(G)]$, by the assumption on $j$.
In the case that $\cof(\g)>\d$,
the Key Lemma \ImprovedKeyLemma\ implies $j\image\g\in M$ since otherwise it would be an
$M$-fresh sequence, and these cannot be added by the gap forcing
$j(\P*\Q)$. So I may assume $\cof(\g)\lte\d$ and consequently
$j\image\g\in M[g]$. Let $C=(P_\k\g)^V$. By 2 this is equal to
$(P_\k\g)^M$. Consequently, $j\image C=\set{j(\s)\st \s\in
C}=\set{j\image\s\st\s\in (P_\k\g)^M}$, and so $j\image C$ is
constructible from $j\image \g$ in $M[g]$.

Let me now argue that $C$ is a $\d$-club in $V[g]$. This part of the
argument is similar to \cite[HamShl].  First, let me argue that $C$ is
unbounded. Any $\s\in (P_\k\g)^{V[g]}$ comes from a name $\dot\s$, and
is therefore covered by $\t=\set{\a\st
\boolval{\check\a\in\dot\s}\not=0}$. It must be that $\t\in C$ since
$\P$ is small. Second, let me show that $C$ is $\d$-closed. Suppose that $D\of
C$ has size less than $\k$ and is $\d$-directed in $V[g]$. I want to show that
$\union D\in C$. It suffices to show that $\union D\in V$. Let $\dot D$
be a name for $D$, and let $D_p=\set{\s\in C\st p\forces \check\s\in
\dot D}$.  Thus, $D=\union_{p\in g}D_p$. There must be some $p\in g$
such that $D_p$ is $\of$-cofinal in $D$ for if not, then for each $p\in
g$ I may choose $\s_p\in D$ such that $D_p$ contains no supersets of
$\s_p$. Since $D$ is $\d$-directed and $\card{g}<\d$ there must be a
$\s\in D$ such that $\s_p\of \s$ for each $p\in g$. But $\s$ must be
forced in $D$ by some condition $p\in g$, so $\s\in D_p$ for some $p\in
g$, contrary to the choice of $\s_p$. Thus, there is some $p\in g$ such
that $D_p$ is $\of$-cofinal in $D$, and consequently $\union D=\union
D_p\in V$, as desired. So $C$ is $\d$-club in $V[g]$.

Therefore $j(C)$ is a $\d$-club in $M[g]$. Also, $j\image C\in M[g]$.
But $j\image C$ is $\d$-directed by $\of$ and has size less than
$j(\k)$.  It follows that $\union (j\image C)\in j(C)$. But $\union
(j\image C)=\union\set{j(\s)\st\s\in C}=\union\set{j\image\s\st\s\in
C}= j\image\g$, and so I conclude that $j\image\g\in j(C)\in M$. Thus,
$j\image\g\in M$, as desired. At the top of the induction, I
conclude $j\image\l\in M$.  So 3 holds.

Now let me prove 4. By 1 and 3, it follows that $j\image \l\in V$.
Suppose now that $j\image\g\in V$, where
$\l\le\g$; I will argue that $j\image\g\plus\in V$ as well (meaning
$(\g\plus)^V$). A simple induction on $n$ will then prove the theorem. By the
Key Lemma \ImprovedKeyLemma\ it suffices to show that every
initial segment of $j\image\g\plus$ is in $V$. So suppose
$\g\le\b<\g\plus$. There must be in $V$ a relation $\triangleleft$ on
$\g$ such that $\b=\ot\<\g,\triangleleft>$. But if $\a$ has order-type
$\z<\b$ with respect to $\triangleleft$, then $j(\a)$ will have
order-type $j(\z)$ with respect to $j(\triangleleft)$. Consequently
$$j\image\b=
 \set{\hbox{order-type of $j(\a)$ wrt $j(\triangleleft)$}\st \a<\g},$$
which is constructible from $j(\triangleleft)$ and $j\image\g$ in $V$,
as desired. Lastly, in the event that $\P$ is $\ltb$-distributive, one
easily obtains the limit stages and I conclude that
$j\restrict\l^{+^{(\a)}}\in V$ for every $\a<\b$.\qed

Please observe that I did not use the hypothesis that $\l<j(\k)$ when proving
1 and 2, so that, for example, those conclusions follow when $j$ is
a huge embedding in $V[G]$. Also, let me point out one important consequence
of the Gap Forcing Theorem. Namely,
if $j:V[G]\to M[j(G)]$ is the $\l$-supercompactness
embedding by a normal fine measure $\eta$ on $(P_\k\l)^{V[G]}$, then
actually $\eta$ concentrates on $(P_\k\l)^V$, and the reason for this is
that $X\in\eta\iff j\image \l\in j(X)$, but $j\image\l\in M$, so
$j\image\l\in P_{j(\k)}j(\l)^M=j(P_k\l^V)$. So every supercompactness
measure added by gap forcing concentrates on a set in the ground model. But
actually, as I will now prove, much much more is true:

\corollary Gap Forcing Corollary. Suppose that $\k$ is supercompact in $V[G]$,
a forcing extension which admits a gap below $\k$. Then $\k$ is supercompact
in $V$, and every supercompactness measure in $V[G]$ extends a measure
in $V$.
\ref\GapForcingCorollary

Let me point out that this corollary is a global claim, in the sense
that the {\it full} supercompactness of $\k$ in $V[G]$ is used to
obtain the {\it full} supercompactness of $\k$ in $V$; and the source
of this phenomenon is identified: every supercompactness measure in
$V[G]$ extends a supercompactness measure in $V$. I will prove the
corollary by proving the following more local version, in which a bit
more than the $\l$-supercompactness of $\k$ in $V[G]$ is used to obtain
the $\l$-supercompactness of $\k$ in $V$. The {\df finite cardinal
successors} of $\l$ are simply the cardinals $\l\plus$, $\l\plusplus$,
$\l\plusplusplus$, and so on, obtained by applying the successor
operation finitely many times.

\corollary Local Improvements. If $V[G]$ admits a gap below $\k$, then:
\points 1. If $\k$ is $2^{\l^\ltk}$-supercompact in $V[G]$ then every
           $\l$-supercompactness measure in $V[G]$ extends a measure in $V$.\cr
        2. If the \GCH\ holds in $V$, or indeed, if only $2^{\l^\ltk}$ is a
           finite cardinal sucessor cardinal of $\l$ in $V$, then again every
           $\l$-supercompactness measure in $V[G]$ extends a measure in $V$.\cr

\ref\LocalImprovements \proof The Gap Forcing Corollary follows from 1,
so let me begin with that. Suppose $\nu$ is a $\l$-supercompactness
measure in $V[G]$. Let $\theta=2^{\l^\ltk}$, and suppose that $\k$ is
$\theta$-supercompact in $V[G]$. Let $\jmu:V[G]\to M_\mu[\jmu(G)]$ be
the ultrapower by a $\theta$-supercompactness measure $\mu$ in $V[G]$.
Since $\card{\nu}=\theta$, it follows that $\jmu\image\nu\in
M_\mu[\jmu(G)]$, and, since this is a subset of $\jmu(\nu)$ of size
less than $\jmu(\k)$, it follows that it has nonempty intersection.
Pick any $s\in\intersect j\image\nu$. It follows that $X\in\nu\iff
s\in\jmu(X)$ for $X\of P_\k\l$. I will use the Gap Forcing Theorem
\GapForcingTheorem\ applied to $\jmu$. Since, by the remarks after the
Gap Forcing Theorem, $\nu$ concentrates on $(P_\k\l)^V$, it follows
that $s\in M_\mu$, and consequently, $s\in V$. Now let
$\<X_\a\st\a<\theta>$ enumerate $P(P_\k\l)^V$ in $V$. Thus,
$\jmu(\<X_\a\st\a<\theta>)\in M_\mu$, and by 1, 2, and 3 of the Gap
Forcing Theorem applied to $\jmu$, I may conclude that
$\<\jmu(X_\a)\st\a<\theta>\in V$. Thus, in $V$ I can simply check
whether $s\in \jmu(X_\a)$ to know whether $X_\a\in\nu$. So
$\nu\intersect V\in V$. It is easy to verify that $\nu\intersect V$ is
a normal fine measure on $P_\k\l$ in $V$. So $\nu$ extends a
supercompactness measure in $V$. So 1 is proved.

To prove 2, I may not assume that $\k$ is anything more than
$\l$-supercompact in $V[G]$, with the measure $\nu$. But I may assume
instead that $\theta=(2^{\l^\ltk})^V$ is a finite successor cardinal to
$\l$ in $V$. Suppose $j:V[G]\to M[j(G)]$ is the ultrapower by $\nu$.
By 4 of the Gap Forcing Theorem, $j\image \theta\in V$, and I can
proceed as before: let $\<X_\a\st \a<\theta>$ enumerate $P(P_\k\l)^V$
in $V$. Since $M\of V$, I know that $j(\<X_\a\st\a<\theta>)\in V$. By
restricting to $j\image\theta$ and collapsing the domain, it follows
that $\<j(X_\a)\st\a<\theta>\in V$. Since $X\in\nu\iff j\image\l\in
j(X)$, I can in $V$ determine which $X_\a$ satisfy $j\image\l\in
j(X_\a)$. Thus, $\nu\intersect V\in V$, and again $\nu$ extends a
measure in $V$, as required.\qed

Just as in the Gap Forcing Theorem \GapForcingTheorem, the \GCH\
assumption in 2 can be relaxed even further when the forcing is
distributive. But the natural question remains, whether in the absence of additional
strength assumed in $V[G]$ or \GCH-type hypotheses in $V$, the
completely local result holds: that is, in a gap forcing extension
$V[G]$, must every cardinal $\k$ which is $\l$-supercompact in $V[G]$
be $\l$-supercompact in $V$? While I don't know the answer, my subsequent
results need only the hypotheses stated here.

\corollary. If $\k$ is supercompact in $V[G]$, a forcing extension which
admits a gap below $\k$, then every $\k$-complete measure in $V[G]$ which
concentrates on a set in $V$ extends a measure in $V$.

\proof Suppose that $\nu$ is a $\k$-complete measure in $V[G]$ which
concentrates on a set $D\in V$. Let $\l=\card{D}$, and let
$j:V[G]\to M[j(G)]$ be a $2^\l$-supercompactness embedding in $V[G]$.
Thus, $j\image 2^\l\in V$ and $M\of V$ by the Gap Forcing Theorem.
It follows that if $\<X_\a\st \a<2^\l>$ enumerates $P(D)$ in $V$, then
$\<j(X_\a)\st \a<2^\l>\in V$. Since $j$ is $2^\l$-supercompact, it follows
that $j\image\nu\in M[j(G)]$, and therefore, since $j(\nu)$ is
$j(\k)$-complete and $2^\l<j(\k)$, that $\intersect j\image\nu\in j(\nu)$. Pick
$s\in \intersect j\image\nu$. Since $s\in j(D)$, it follows that
$s\in M$, and hence also $s\in V$. Also, $X\in\nu\iff s\in j(X)$. Therefore,
I can tell in $V$ whether $X_\a\in\nu$ by checking whether $s\in j(X_\a)$.
So $\nu\intersect V\in V$, as desired.\qed

\corollary. Suppose that $2^\k$ is a finite cardinal successor of $\k$,
that $V[G]$ admits a gap below $\k$, and that $\k$ is measurable in $V[G]$.
Then every measure on $\k$ in $V[G]$ extends a measure in $V$.

\proof For normal measures on $\k$, this is a special case of the Local
Improvements \LocalImprovements\ of the Gap Forcing Corollary. But
nearly the same argument works for any measure on $\k$.  Suppose $\nu$
is a measure on $\k$ in $V[G]$, with embedding $j:V[G]\to M[j(G)]$. By
the Gap Forcing Theorem \GapForcingTheorem\ it follows that $M\of V$
and $j\image 2^\k\in V$, since $2^\k=\k^{{+}^{(n)}}$ for some finite $n$. It
follows that $\<j(X_\a)\st\a<2^\k>\in V$ where $\<X_\a\st\a<2^\k>$
enumerates $P(\k)$ in $V$. Since $X\in\nu\iff [\id]_\nu\in j(X)$, I can
decide in $V$ whether $X_\a\in\nu$. So $\nu\intersect V\in V$, and thus
$\nu$ extends a measure on $\k$ in $V$.\qed

Let us say that a poset admits a {\df very low} gap if it admits a gap
at or below the least inaccessible cardinal. (This could be modified,
without affecting any of the arguments below, to the least Mahlo
cardinal, the least weakly compact cardinal, or indeed anything that is
strictly below what might become the least measurable cardinal in a
forcing extension; e.g. the least weakly compact limit of weakly
compact cardinals.)

\corollary No Turn-On Corollary. If $\P$ admits a very low gap, then it creates
no supercompact cardinals. And, if the \GCH\ holds in $V$, neither does it
increase the degree of supercompactness of any cardinal. In particular, it
creates no new measurable cardinals.

\proof This follows immediately from the Local Improvements
\LocalImprovements\ of the Gap Forcing Corollary, since $\P$ admits a
gap below all the cardinals in question. And I don't actually need the
full \GCH; rather, to show every $\l$-supercompactness measure in
$V[G]$ extends a measure in $V$, I need only know that $2^{\l^\ltk}$ is
a finite successor cardinal of $\l$, and even less than this if $\P$ is
distributive.\qed

The Gap Forcing Theorem \GapForcingTheorem\ and its corollaries point
at the following conjecture, which, if true, would explain and unify
its conclusions.

\conjecture Gap Forcing Conjecture. After forcing with a gap below $\k$, every
supercompactness embedding with critical point $\k$ lifts an embedding
from the ground model.

\noindent The conjecture asserts that if $j:V[G]\to M[j(G)]$ is a
supercompactness embedding with critical point $\k$, then $j\restrict
V$ is a definable class in $V$. The conjecture does not assert that
$j\restrict V$ is a supercompactness embedding in $V$; I have already
pointed out, in my remarks preceding the Gap Forcing Theorem,
how that can fail. An affirmative answer to the following question would
strongly unify the consequences of the Gap Forcing Theorem:

\question. After forcing which admits a gap below $\k$, does every
ultrapower embedding by a $\k$-complete measure on any set lift an
embedding from the ground model?

A even stronger version of this question asks: after forcing which admits a
gap below $\k$, does every embedding with critical point $\k$ lift an
embedding from the ground model? But this is too strong. A
counterexample is provided by the situation in which there is a
measurable cardinal $\k$ with two distinct normal measures $\mu$ and
$\nu$. Add a Cohen real $x$ and then add a Cohen subset to $\omega_2$.
This admits a gap at $\omega_1$, and the measures $\mu$ and $\nu$ lift
and extend uniquely to measures $\bar\mu$ and $\bar\nu$, respectively,
in $V[G]$. The embedding $j:V[G]\to M[j(G)]$ obtained by iterating the
measures $\bar\mu$ and $\bar\nu$ in the order specified by $x$ cannot
be the lift of an embedding definable in $V$, because from
$j\restrict V$ one can recover the real $x$. Thus, after forcing which
admits a gap below $\k$, it needn't be that every embedding with
critical point $\k$ lifts an embedding from the ground model.

\section Partial Laver Preparations

Laver's crucial contribution in \cite[Lav78] was the idea of what is now
called a {\df Laver function}. For a supercompact cardinal $\k$, this
is a function $\ell:\k\to V_\k$ such that for any $\l\geq\k$ and any
$x\in H(\l\plus)$ there is a $\l$-supercompact embedding $j:V\to M$
such that $j(\ell)(\k)=x$. The function $\ell$ can be defined
inductively:  for a measurable non-supercompact $\g$, let $\l$ be least
such that for some $x\in H(\l\plus)$ there is no $\l$-supercompact
embedding $j:V\to M$ with critical point $\g$ such that
$j(\ell\restrict\g)(\g)=x$, select a witness $x$, and let
$\ell(\g)=x$.  Since this definition is local, in the sense that
$\ell(\g)$ depends only on choices made concerning $\ell\restrict\g$,
and not on the ultimate length of the function $\ell$ to be defined, it
actually gives a class function $\ell$, called the {\df universal Laver
function}, whose initial segments work as a Laver function
simultaneously for every supercompact cardinal (see \cite[KimMag],
\cite[Apt96] for elaboration).  There are several other simplifying
assumptions to be made about a Laver function $\ell$, and in this paper
I will take them to be part of the definition of what it means to be a
Laver function. First of all, I may assume that $\dom(\ell)$ consists
entirely of measurable non-supercompact cardinals. Furthermore, I may
assume that if $\g\in\dom(\ell)$ then $\g$ is closed under $\ell$ in
the sense that $\ell\image\g\of V_\g$.  Third, I may assume, as I
point out in \cite[Ham94a], that for every $\l$ and every $x\in
H(\l\plus)$ there is a $\l$-supercompact embedding $j:V\to M$ such that
$j(\ell)(\k)=x$ and $\dom(j(\ell))\intersect(\k,\l]=\emptyset$. Thus, a
Laver function has long gaps in its domain.

The {\df Laver preparation} is the forcing iteration defined by Laver
\cite[Lav78] from a Laver function $\ell$. It has reverse Easton support,
so that direct limits are taken at inaccessible stages, and inverse
limits otherwise. The stage $\g$ forcing is exactly $\ell(\g)$,
provided, as perhaps seems unlikely, but as actually occurs on a
stationary set of $\g$, that this is the $\P_\g$ name of a
$\ltg$-directed closed poset in $V^{\P_\g}$.  If $\ell(\g)$ is not such
a name, then the stage $\g$ forcing is trivial.  Laver proved that
performing the Laver preparation makes a supercompact cardinal $\k$
{\df indestructible}, so that any further $\ltk$-directed closed
forcing will preserve the supercompactness of $\k$. The {\df universal}
Laver preparation is simply the Laver preparation obtained by using a
universal Laver function, and it has the effect of making every
supercompact cardinal indestructible \cite[KimMag], \cite[Apt96].

In what may be the obvious choice of poset, given that I want to make a
supercompact cardinal {\it partially} indestructible, I define that a
poset $\P$ is a {\df partial} Laver preparation of $\k$ iff there is a
Laver function $\ell$ on $\k$ and a set $A\of\dom(\ell)$, the set of
{\df allowed} stages, such that $\P$ is the reverse Easton
$\k$-iteration which at stage $\g$ forces with the poset
$\dot\Q_\g=\ell(\g)$, provided, as usual, that $\ell(\g)$ is
the $\P_\g$-name of a poset which is
$\ltg$-directed closed in $V^{\P_\g}$, but also, that $\g$ is an
allowed stage, i.e. that $\g\in A$. It is clear that any partial Laver
preparation admits a gap between any two allowed stages.

\lemma Partial Laver Preparation Lemma. After a partial Laver
preparation of a supercompact cardinal $\k$, it remains supercompact.
\ref\PartialLaverPreparationLemma

\proof The usual Laver \cite[Lav78] argument adapts to this
circumstance.  Suppose $A\of\k$ is the set of allowed stages, and
assume that $G\of\P$ is $V$-generic.  If $A$ is bounded below $\k$,
then the forcing is small, and so $\k$ remains supercompact. Otherwise,
the generic $G$ has size $\k$. Fix any $\l$, and $\theta\muchgt\l$, and
pick a $\theta$-supercompact embedding $j:V\to M$ such that
$j(\ell)(\k)$ is not the name of a poset, and $\dom(j(\ell))\intersect
(\k,\theta]=\emptyset$. Thus, since there can be no allowed stages from
$\k$ to $\theta$, I can factor $j(\P)$ as $\P*\Ptail$, where $\Ptail$
is ${\lte}\theta$-closed in $M[G]$. Thus, it is also
${\lte}\theta$-closed in $V[G]$. Let $\Gtail\of\Ptail$ be
$V[G]$-generic, and lift the embedding in $V[G][\Gtail]$ to
$j:V[G]\to M[j(G)]$ where $j(G)=G*\Gtail$. Now, use $j\image\l$ as a
seed to define a measure $\mu$ on $P_\k\l$ according to the rule
$X\in\mu\iff j\image\l\in j(X)$. It is straightforward to check that $\mu$ is a
normal fine measure on $P_\k\l$. Furthermore, $\mu$ must be in $V[G]$,
since it cannot have been added by the closed forcing $\Gtail$. Thus,
$\k$ remains $\l$-supercompact in $V[G]$.\qed

\lemma No Turn-on Lemma. Except possibly for the cardinal $\g$ which is
the first nontrivial stage of forcing, a partial Laver preparation adds
no supercompact cardinals. If the \GCH\ holds in the ground model,
then, except again for $\g$, neither does it increase the degree of
supercompactness of any cardinal. (By adding a Cohen real in front of
the preparation, or any other small forcing, these provisos about
$\g$ can be removed.)

\proof This follows immediately from the Gap Forcing Corollary
\GapForcingCorollary\ and the Local Improvements \LocalImprovements\ in
the previous section, since the partial Laver preparation admits a gap
just after its first stage.  By adding a Cohen real in front, or indeed any
small enough forcing, a very low gap is introduced, and the conclusion
applies even at $\g$.\qed

Without the forcing in front of the preparation, it is possible that a
partial Laver preparation could increase the degree of supercompactness
of its very first stage $\g$. This would occur, for example, if $\g$
had been previously supercompact, and the current model had been
obtained by forcing with reverse Easton support to add subsets to the
measurable elements of a certain club of cardinals below $\g$. It can be
arranged that $\g$ would be measurable but not $2^\g$-supercompact
after this forcing, but then become suddenly supercompact again after
the first stage of the subsequent partial preparation.

\section Jumping High

The concept of a high jump function will prove useful
later on in making the tail of an iteration sufficiently closed.
A {\df high jump} function for a (partially) supercompact
cardinal $\k$ is a function $h\from\k\to\k$ such
that $j(h)(\k)>\l$ whenever $j$ is a $\l$-supercompact embedding on $\k$.

\lemma High Jump Lemma. Suppose $\k$ is supercompact but no
normal measure on $\k$ concentrates on the supercompact cardinals below
$\k$.  Then there is a high jump function for $\k$. In particular, the
least supercompact cardinal has a high jump function.

\proof If $\b<\k$ is not supercompact, let $h(\b)=2^{\z^{<\b}}$, where $\z$ is
least such that $\b$ is not $\z$-supercompact. Since $\k$ is supercompact,
it follows that $h(\b)<\k$ and so $h\from\k\to\k$. Now suppose $j:V\to M$ is a
$\l$-supercompact embedding. Necessarily $\k$ is not supercompact in
$M$, since the induced normal measure does not concentrate on the
supercompact cardinals below $\k$. Thus $\k$ fails in $M$ to be
$\z$-supercompact for some minimal $\z$.
If $2^{\z^\ltk}\le\l$ in $M$, then I may code a $\z$-supercompactness measure
on $P_\k\z$ in $V$ with a subset of
$\l$. Since $P(\l)^M=P(\l)^V$, this set, and hence also the measure, must be in
$M$, a contradiction. Therefore, $j(h)(\k)=(2^{\z^\ltk})^M>\l$, as desired.\qed

\lemma Another High Jump Lemma. If there are fewer than $\k$ many measurable cardinals
above the supercompact cardinal $\k$, then $\k$ has a high-jump function.

\proof Suppose there are exactly $\a$ many measurable cardinals above $\k$,
and $\a<\k$. For $\g<\k$,
let $h(\g)$ be the $(\a+1)^\th$ measurable cardinal above $\g$.
Thus, $h:\k\to \k$, and
if $j:V\to M$ is a $\l$-supercompact embedding, then $j(h)(\k)$ is the
$(\a+1)^\th$ measurable cardinal above $\k$ in $M$. Up to $\l$, however, if a
cardinal is measurable in $M$ then it is measurable in $V$, so the
$(\a+1)^\th$ measurable cardinal in $M$ is above $\l$. That is,
$j(h)(\k)>\l$, as desired.\qed

This technique can be pushed much harder. For example, if there are
$\k$ many measurable cardinals above $\k$, one uses the function in which
$h(\g)$ is the $\g^\th$ measurable cardinal above $\g$. If there are
$\k\plus$ many, let $h(\g)$ be the ${\g\plus}^\th$ measurable above
$\g$. These ideas lead naturally to the ideas of \S 5. Also, though,
in a different sort of generalization, $\k$ need not be {\it fully}
supercompact; if, for example, $\k$ is partially supercompact (but still
a limit of measurables), then the function $h$ defined in the proof will
be a high-jump function for embeddings up to the degree of supercompactness
of $\k$.

The next theorem shows how the existence of a high jump function is
robust.

\lemma High Jump Preservation Lemma. Suppose that $h$ is a high jump
function for $\k$ in $V$. Then it remains so in any forcing extension
$V[G]$ in which every supercompactness measure extends a measure in $V$.
In particular, if $V[G]$ admits a gap below $\k$ and either the
\GCH\ holds in $V$ or $\k$ is sufficiently supercompact in $V[G]$, then
$h$ is a high-jump function for $\k$ in $V[G]$.
\ref\HighJumpPreservationLemma

\proof Suppose that $\nu$ is a $\l$-supercompactness measure in $V[G]$,
with the corresponding embedding $j:V[G]\to M[j(G)]$, and
$\nu$ extends a
$\l$-supercompactness measure $\mu$ in $V$.  Let
$X=\set{j(f)(j\image\l)\st f\in V}\elesub M$ be the seed hull of
$j\image\l$ via $j\restrict V$. Since $X$ is isomorphic to $M_\mu$, the
ultrapower of $V$ by $\mu$, by the map $\varphi:j(f)(j\image\l)\mapsto
[f]_\mu$, it follows that $\jmu=\pi\compose j$, where $\pi$ is the
collapse of $X$ (See \cite[Ham97b] for elaboration on this seed hull
factor method). Consequently, since $\l<\jmu(h)(\k)$, it follows that
$\l<\pi(j(h))(\k)=\pi(j(h)(\k))\leq j(h)(\k)$, as desired.

If either the \GCH\ holds in $V$ or $\k$ is sufficiently supercompact
in $V[G]$, then the Gap Forcing Corollary \GapForcingCorollary\ and
the Local Improvements \LocalImprovements\ yield
the necessary hypothesis that every supercompactness measure in
$V[G]$ extends a measure in $V$. For $h$ to work with $\l$-supercompactness
embeddings in $V[G]$, one needs to know either that $2^{\l^{<\k}}$ is a
finite cardinal successor to $\l$ in $V$ (less if the forcing is
distributive) or that $\k$ is
$2^{\l^{<\k}}$-supercompact in $V[G]$.\qed

One might suspect that there is a high jump function for any
supercompact cardinal; but the following theorem should temper that
tendency.

\theorem Almost Huge High Jump Theorem. If $\k$ is almost huge, then
there is no high jump function for $\k$.

\proof Suppose $j:V\to M$ is an almost huge embedding, so that $M^{<j(\k)}\of M$. It
follows that $j(\k)$ is regular, and so the set $\set{j(h)(\k)\st h\in \k^\k}$,
which has size $2^\k$, is bounded by some $\l<j(\k)$.
Let $\mu$ be the supercompactness embedding germinated
by the seed $j\image\l$, and let $\pi$ be the collapse of the seed hull
$X=\set{j(f)(j\image\l)\st f\in V}\elesub M$. It follows that
$\jmu=\pi\compose j$, and thus, for any function $h\in\k^\k$ we have
$\jmu(h)(\k)=\pi(j(h))(\k)=\pi(j(h)(\k))<\pi(\l)\leq\l$. So there can be no
high jump function that works with the measure $\mu$.\qed

Nevertheless, I will show how to add high jump functions for every
supercompact cardinal, and moreover, to do this in a way which preserves
all supercompact cardinals. By the previous theorem, all almost huge
cardinals will of necessity be destroyed.
In this argument, I will use the concept of coherent
clubs, which were first introduced by Hugh Woodin \cite[W] in his construction
to obtain a model of a
supercompact cardinal whose weak compactness is easily destroyed. They
later returned with a vengeance---huge infinities of them piling up all
around---as the central technique, and the central complication, of the
Fragile Measurability Theorem \cite[Ham94a]. They appear also in the
epilogue of this paper, when I use them to separate the notions of
fragility and superdestructibility.

\lemma Coherent Club Lemma. While preserving all supercompact
cardinals, and in fact making them indestructible, one can add by
forcing a system of clubs $C_\g\of\g$ for inaccessible cardinals $\g$,
each disjoint from the supercompact cardinals. Furthermore,
the clubs can be made to cohere in the sense that if $\d$ is an
inaccessible cluster point of $C_\g$, then $C_\g\intersect\d=C_\d$.
The forcing does not create any supercompact cardinals, and if the
\GCH\ holds in the ground model, neither does it increase the
degree of supercompactness of
any cardinal; in fact, every new supercompactness measure extends an old
measure.

\proof I will interweave the universal Laver preparation with the
forcing to add a system of coherent clubs avoiding the supercompact
cardinals.  Specifically, let $\ell$ be a universal Laver function, and
let $\P$ be the following reverse Easton class iteration:
at inaccessible stages $\g$, the stage $\g$ forcing has two parts. First
is the coherent club forcing $\Q_\g$, whose conditions are closed bounded
subsets $C\of\g$, ordered by end-extension, such that $C$ contains no
supercompact cardinals of $V$, and if $\d$ is an inaccessible cluster point
of $C$, then $C\intersect\d=C_\d$, the club added earlier at stage $\d$.
This forcing has open dense sets as closed as you like up to $\g$: the set
of conditions mentioning a point above $\l<\g$ is $\lte\l$-directed closed.
The second part of the stage $\g$ forcing, $\R_\g$, is simply the forcing
given by the
Laver function $\ell(\g)$, if this is the name of $\ltg$-directed
closed poset in $V^{\P_\g*\Q_\g}$ (if not, then $\R_\g$ is trivial). Please
observe that this iteration admits a very low gap.

Suppose now that $G$, a proper class, is $V$-generic for the forcing
$\P$. Let me show that every supercompact cardinal of $V$ is preserved
to $V[G]$. Suppose that $\k$ is supercompact in $V$. I will show that
in fact $\k$ becomes indestructible in $V[G]$. So suppose $H\of\Q$ is
$V[G]$-generic for the $\ltk$-directed closed forcing $\Q$. Fix $\l$,
and pick $\theta\muchgt\l,\card{\Q}$. It suffices to show that $\k$ remains
$\l$-supercompact after forcing with $\P_\theta*\Q$. Factor $\P_\theta$
as $\P_\k*\Q_\k*\P_{\k,\theta}$, where $\Q_\k$ is the stage $\k$
coherent club forcing, and $\P_{\k,\theta}$ is the rest of the forcing
up to stage $\theta$, beginning with $\R_\k$.  Thus, $\P_{\k,\theta}$ has a
$\ltk$-directed closed dense set
in $V^{\P_\k*\Q_\k}$. I may replace $\P_{\k,\theta}$ with this dense set and
assume that $\P_{\k,\theta}$ itself is $\ltk$-directed closed.
Fix $j:V\to M$ a $\theta$-supercompact embedding
in $V$ such that $j(\ell)(\k)=\P_{\k,\theta}*\Q$ and
$\dom(j(\ell))\intersect(\k,\theta]=\emptyset$. In particular, $\k$ is
not supercompact in $M$, since $\dom(\ell)$ contains no supercompact
cardinals. Observe that
$j(\P_\k)=\P_\k*\Q_\k*(\P_{\k,\theta}*\Q)*\Ptail=P_\theta*\Q*\Ptail$
where $\Ptail$ is ${\lt}\theta$-closed in $M^{\P_\theta*\Q}$, and hence
also in $V[G][H]$. Let $\Gtail$ be $V[G][H]$-generic for $\Ptail$, and
then lift to $j:V[G_\k]\to M[j(G_\k)]$ where
$j(G_\k)=G_\k*C_\k*(G_{\k,\theta}*H)*\Gtail$. That is,
$j(G_\k)=G_\theta*H*\Gtail$. Now let $\bar C_\k=C_\k\union\{\k\}$.
This is a condition in $j(\Q_\k)$ since, first, $\k$ is not
supercompact in $M$, second, the reflection $\bar
C_\k\intersect\k=C_\k$ is the generic used at stage $\k$ in $j(G_\k)$, and,
third, the reflection property holds below $\k$ since $C_\k$ was generic.
To use suggestive notation, let $j(C_\k)$ be $V[G][H][\Gtail]$-generic for
$j(\Q_\k)$ below the master condition $\bar C_\k$. This forcing has a
dense set which is $\lte\theta$-directed closed. Now lift the embedding, in
$V[G][H][\Gtail][j(C_\k)]$, to $j:V[G_\k][C_\k]\to
M[j(G_\k)][j(C_\k)]$.  Observe that $j\image G_{\k,\theta}\in
M[j(G_\k)][j(C_\k)]$, and $j(\P_{\k,\theta})$ is ${\lt}j(\k)$-directed
closed in that model. Thus, I can find a master condition below
$j\image G_{\k,\theta}$, force below it, and lift to
$j:V[G_\k][C_\k][G_{\k,\theta}]\to
M[j(G_\k)][j(C_\k)][j(G_{\k,\theta}]]$. That is, $j:V[G_\theta]\to
M[j(G_\theta)]$. Finally, use $j\image H$ as a master condition, add a
generic $j(H)$,  and lift to $j:V[G_\theta][H]\to
M[j(G_\theta)][j(H)]$. This lift is defined in
$V[G_\theta][H][\Gtail][j(C_\k)][j(G_{\k,\theta})][j(H)]$. Using $j\image\l$
as a seed, I germinate a normal fine measure $\mu$ on
$(P_\k\l)^{V[G_\theta][H]}$. Since the extra tail forcing and master
condition forcing was ${\lt}\theta$-closed, it could not have added
$\mu$, and consequently $\mu$ lives in $V[G_\theta][H]$, as desired. So
every supercompact cardinal in $V$ remains supercompact---becoming in fact
indestructible---in $V[G]$.

Since $\P$ admits a very low gap, the No Turn-On Lemma tells us that it cannot
create any supercompact cardinals, and, if the \GCH\ holds, neither can it
increase the degree of supercompactness of any cardinal; every new supercompactness
measure extends an old measure. Thus, since the clubs
which I added are disjoint from the supercompact cardinals of $V$,
they are also disjoint from the supercompact cardinals of $V[G]$. Also, since
the clubs are built from initial segments with the coherence property,
the clubs themselves also have the coherence property.\qed

\theorem Universal High Jump Theorem. While preserving all supercompact
cardinals and creating no new supercompact cardinals, one can, via forcing
with a very low gap, add high jump functions for every supercompact
cardinal.

\proof If every inaccessible cardinal has a club subset disjoint from the
smaller supercompact cardinals, then
no normal measure can concentrate on a set of supercompact cardinals. So
the corollary follows from the previous lemma and the High Jump Lemma.
In fact, after adding the clubs, the proof of the High Jump Lemma produces a single
class function
$h\from\ORD\to\ORD$ whose restriction $h\restrict\k$ to any supercompact
cardinal $\k$ yields a high jump function for $\k$.\qed

\theorem High Jump Theorem. While preserving all supercompact cardinals and
creating no new supercompact cardinals, via forcing with a very low
gap, one can add high jump functions for every supercompact cardinal up to and
including $\k$, without collapsing cardinals above $\k$. If the \GCH\ holds,
then it can also be arranged to add the functions without collapsing cardinals
or cofinalities at all.
\ref\HighJumpTheorem

\proof For the first part, just perform the coherent club forcing up to and including
the stage $\k$ coherent club forcing.
This forcing has size $\k$, and hence preserves all cardinals above $\k$. For the
second part, when the \GCH\ holds, one can perform a modified coherent
club forcing in which the stage $\g$ forcing is allowed only when, in addition,
it preserves all cardinals, cofinalities, and the \GCH. I will make a
similar argument---giving all the details---in Theorem 4.3.\qed

\section Exact Preservation Theorems

Let me now set off into the unknown continent between the extremes of
indestructibility and superdestructibility. My initial explorations
will reveal the rich structure to be found there:
exact preservation theorems. In these theorems, borrowing from
both indestructibility and superdestructibility, I will precisely control
the class of $\ltk$-directed closed posets which preserve the
supercompactness of $\k$, and obtain models where {\it exactly} a
certain class of posets preserves the supercompactness of a given
supercompact cardinal $\k$. I will lean on Laver's methods to show that
a poset preserves supercompactness, and on the Gap Forcing Theorem
\GapForcingTheorem\ and
its corollaries to show that a poset destroys supercompactness.

Let me begin with two warm-up theorems, in which I show
that the notions of superdestructibility at $\k$
and at $\k\plus$ are orthogonal. Recall from \cite[Ham97b] that a
supercompact cardinal $\k$ is {\df superdestructible at $\theta$} when
any $\ltk$-closed forcing which adds a subset to $\theta$ destroys the
$\theta$-supercompactness of $\k$.

\theorem Exact Preservation Theorem. Assume that $\k$ is supercompact
in $V$. Then there is a forcing extension $V[G]$ in which $\k$ remains
supercompact and becomes superdestructible at $\k\plus$ but not at
$\k$.

\proof I will obtain a model where $\k$ is supercompact and any
$\ltk$-closed poset which adds a subset to $\k\plus$ destroys the
$\k\plus$-supercompactness of $\k$, but where the measurability of $\k$
is preserved by any $\ltk$-directed closed poset which preserves
$\k\plus$ and $2^\k=\k\plus$.

To begin I may assume, by forcing if necessary,
that $2^\k=\k\plus$ in $V$. After this, let $\P$ be the
partial Laver preparation of $\k$ in which the stage $\g$ forcing is allowed
only when it destroys the measurability of $\g$, but preserves
$\g\plus$ and $2^\g=\g\plus$, if indeed this held in $V$.
Suppose that $G\of\P$ is $V$-generic. By
the Partial Laver Preparation Lemma, I know that $\k$ remains
supercompact in $V[G]$.

Let me now prove that $\k$ is superdestructible at $\k\plus$ in
$V[G]$.  Suppose that $H\of\Q$ is $V[G]$-generic, where $\Q$ is a
$\ltk$-closed forcing notion which adds a new subset $B\of\k\plus$, but that
$\k$ remains $\k\plus$-supercompact in $V[G][H]$. (Here I mean
$\k\plus$ to denote $(\k\plus)^V=(\k\plus)^{V[G]}$.) Thus, there is a
$\k\plus$-supercompact embedding $j:V[G][H]\to M[j(G)][j(H)]$. By
the Gap Forcing Theorem \GapForcingTheorem,
I know $P(\k\plus)^M=P(\k\plus)^V$, and thus $2^\k=\k\plus$ in $M$. If
$\k\plus$ is collapsed by $H$, then it must also be collapsed by
$j(G)$, which is impossible since the only stage which could do this is
the stage $\k$ forcing, and that stage is not allowed if it collapses
$\k\plus$. Thus, $H$ must not collapse $\k\plus$.  If the stage $\k$
forcing is not allowed, then $j(G)=G*\Gtail$, where $\Gtail$ is
${\lte}{\k\plus}$-closed, and consequently $B\in M[G]\of V[G]$, a
contradiction. If the stage $\k$ forcing is allowed, then
$j(G)=G*g*\Gtail$, where $g$ is the stage $\k$ forcing and $\Gtail$ is
again ${\lte}{\k\plus}$-closed. Since $\k$ was allowed I know that $\k$
is not measurable in $M[G][g]$. But
$P(\k\plus)^{M[G][g]}=P(\k\plus)^{M[j(G)][j(H)]}=P(\k\plus)^{V[G][H]}$,
and so by coding a measure on $\k$ from $V[G][H]$ into a subset of $\k\plus$
I conclude that it lies in $M[G][g]$, and so $\k$ is measurable
there after all, a contradiction.

Finally, I will show that in $V[G]$ the measurability of $\k$ is
preserved by any $\ltk$-directed closed forcing which preserves
$\k\plus$ and $2^\k=\k\plus$. Suppose $H\of\Q$ is generic for such a
poset, but that $\k$ is not measurable in $V[G][H]$. Fix a large
$\theta$ and a $\theta$-supercompact embedding $j:V\to M$ such that
$j(\ell)(\k)=\dot\Q$ and
$\dom(j(\ell))\intersect(\k,\theta]=\emptyset$. Since $\k$ is not
measurable in $V[G][H]$, this is also true in $M[G][H]$, and,
similarly, $H$ preserves $\k\plus$ and $2^\k=\k\plus$ over $M[G]$. Thus,
$\k$ is an allowed stage in the $j(\P)$ forcing, and so $j(\P)$ factors
as $\P*\Q*\Ptail$ where $\Ptail$ is ${\lte}\theta$-closed. Force over the
tail and lift the embedding to $j:V[G]\to M[j(G)]$ where $j(G)=G*H*\Gtail$ in
$V[G][H][\Gtail]$. Observe that $j\image H\in M[j(G)]$, and so by the
directed closure of $j(\Q)$, there is a master condition $p\in j(\Q)$ which is
below every element of $j\image H$. Let $j(H)$ be generic
below $p$, and lift the embedding to $j:V[G][H]\to M[j(G)][j(H)]$ in
$V[G][H][\Gtail][j(H)]$. Using $\k$ as a seed, I obtain a normal measure
$\mu$ on $P(\k)^{V[G][H]}$. By closure considerations I
know $\mu$ must be in $V[G][H]$, and so $\k$ is measurable there.\qed

Next, for the second warm-up theorem,
I will prove the opposite combination. Define that a supercompact
cardinal $\k$ is {\df indestructible above $\theta$} iff any
$\ltk$-directed closed forcing which adds no subsets to $\theta$
preserves the supercompactness of $\k$.

\theorem Exact Preservation Theorem. Assume $\k$ is supercompact in
$V$. Then there is a forcing extension $V[G]$ in which $\k$ remains
supercompact and becomes superdestructible at $\k$ but not at
$\k\plus$. Indeed, $\k$ can be made simultaneously superdestructible at
$\k$ and indestructible above $\k$.

\proof  Here I will perform the partial Laver preparation of $\k$
in which the
stage $\g$ forcing is allowed provided that it adds no new subsets to
$\g$. Let me emphasize here for clarity that it is part of the definition of
a partial Laver preparation that, in addition, the forcing $\Q_\g$ which
the Laver function hands to us at stage $\g$ must be
$\ltg$-directed closed in $V^{\P_\g}$.
Suppose that $G\of \P$ is $V$-generic for this forcing. By the
Partial Laver Preparation Lemma, it follows that $\k$ remains
supercompact in $V[G]$.

If $\k$ is not superdestructible at $\k$ in $V[G]$, then there must be
some $\ltk$-closed forcing $\Q$, and a $V[G]$-generic $H\of\Q$, adding
a new subset $B\of\k$, such that $\k$ remains measurable in $V[G][H]$,
with the corresponding embedding $j:V[G][H]\to M[j(G)][j(H)]$. By
closure considerations, $B\in M[G]$, since the stage $\k$ forcing
cannot have added a subset to $\k$. But the Gap Forcing
Theorem \GapForcingTheorem\ tells us that $M[G]\of V[G]$, and so
$B\in V[G]$, contradicting the fact that it was newly added by $H$.
So $\k$ is superdestructible in $V[G]$, as desired.

To show that $\k$ becomes indestructible above $\k$ in $V[G]$, I will
employ what I will later refer to as the `usual' lifting argument: suppose
$H\of\Q$ is $V[G]$-generic for $\ltk$-directed closed forcing $\Q$
which adds no new subsets to $\k$. Fix any $\l$, and pick
$\theta\muchgt\l,\card{\Q}$. Select a $\theta$-supercompact embedding
$j:V\to M$ such
that $j(\ell)(\k)=\dot\Q$ and $\dom(j(\ell))\intersect
(\k,\theta]=\emptyset$. The stage $\k$ is allowed since $\Q$ adds no subsets
to $\k$. Thus, $j(\P)=\P*\Q*\Ptail$, where $\Ptail$ is ${\lte}\theta$-closed.
Force to add a generic $\Gtail\of\Ptail$, and lift the embedding to
$j:V[G]\to M[j(G)]$, where $j(G)=G*H*\Gtail$. The lift is defined in
$V[G][H][\Gtail]$. Now $j(\Q)$ is ${\lt}j(\k)$-directed closed, and
$j\image H\in M[j(G)]$, so I can find a master condition below $j\image H$
in $j(\Q)$, and force below it to add a generic $j(H)\of j(\Q)$. This
gives a lift embedding $j:V[G][H]\to M[j(G)][j(H)]$, defined in
$V[G][H][\Gtail][j(H)]$. Use $j\image\l$ as a seed to germinate a normal
fine measure
$\mu$ on $P_\k\l^{V[G][H]}$ as follows: $X\in \mu\iff j\image\l\in j(X)$.
Since the forcing $\Gtail*j(H)$ was ${\lte}\theta$-closed, it could not
have added $\mu$, and consequently $\mu\in V[G][H]$. Thus, $\k$ is
$\l$-supercompact there, as desired.\qed

I can in addition exhibit models of the third and fourth possibilites:
after small forcing, the main result of \cite[Ham97b] shows that a
supercompact cardinal is superdestructible at both $\k$ and $\k\plus$;
and if $\k$ is indestructible, then $\k$ is superdestructible at
neither $\k$ nor $\k\plus$. Let me now continue with additional
Exact Preservation Theorems.

\theorem Exact Preservation Theorem. Assume $\k$ is supercompact and the
\SCH\ holds in $V$. Then
there is a forcing extension $V[G]$, obtained without collapsing
cardinals or cofinalities, in which $\k$ remains supercompact
and indeed over
$V[G]$ the supercompactness of $\k$ is preserved by exactly those
$\ltk$-directed closed posets which collapse neither cardinals nor
cofinalities.
\ref\EPCard

\proof Let $\P$ be the partial Laver preparation of $\k$
in which stage $\g$ is
allowed provided that it collapses neither cardinals nor cofinalities,
and suppose that $G\of\P$ is $V$-generic. By the Partial Laver
Preparation Lemma, $\k$ remains supercompact in $V[G]$.

Let me argue that the supercompactness of $\k$ is preserved by any
$\ltk$-directed closed forcing which collapses neither cardinals nor
cofinalities. Suppose $\Q$ is such forcing, and $H\of\Q$ is
$V[G]$-generic. Fix $\l$, pick $\theta\geq2^\l,\card{\Q}$, and select as usual
$j:V\to M$ a $\theta$-supercompact embedding such that
$j(\ell)(\k)=\Qdot$ and $\dom(j(\ell))\intersect
(\k,\theta]=\emptyset$.  Since $H$ preserves cardinals and cofinalities
over $V[G]$, it follows by the largeness of $\theta$ that this is also
true over $M[G]$, and so the stage $\k$ forcing is allowed. By the usual
lifing argument, lift to
$j:V[G][H]\to M[j(G)][j(H)]$ and conclude that $\k$ remains
$\l$-supercompact in $V[G][H]$.

Conversely, I will also argue that the supercompactness of $\k$ is
destroyed by any $\ltk$-closed forcing $\Q$ which collapses cardinals
or cofinalities. Necessarily, $\Q$ must collapse the cofinality of some regular
cardinal. Suppose $H\of\Q$ is $V[G]$-generic and, for some
regular cardinal $\l$ in $V[G]$, I have $\cof(\l)=\d<\l$ in $V[G][H]$,
but that $\k$ remains $\d$-supercompact in $V[G][H]$. Let $j:V[G][H]\to
M[j(G)][j(H)]$ be the witness embedding. By the closure of the embedding
it follows that $\cof(\l)^{M[j(G)][j(H)]}=\d$. By the Gap Forcing Theorem
\GapForcingTheorem\ I
know that $M$ and $V$ have the same $\d$-sequences of ordinals, and
consequently, $\cof(\l)^M>\d$. But $j(G)$ cannot have collapsed any
cardinals or cofinalities over $M$, since, as I will prove in the next
paragraph, $G$ did not over $V$, and
$j(H)$ is ${<}j(\k)$-closed, so it cannot have added the $\d$-sequence
either, a contradiction.

I must now prove that neither cardinals nor cofinalities were
collapsed between $V$ and $V[G]$. Here I will use the \SCH\ assumption,
but it is not onerous. The \SCH\ follows of course from the \GCH, which
one can easily force while preserving supercompactness. Also, though,
Solovay \cite[Sol74] proved that the \SCH\ holds automatically above
any supercompact cardinal $\k$, and by reflection it must hold
unboundedly often below $\k$.  The content of the hypothesis is merely
that the \SCH\ holds at the remainder of the singular cardinals below
$\k$. Certainly if $\k$ is not the least supercompact cardinal, then I
could have omitted the \SCH\ assumption entirely, by starting the partial Laver
preparation beyond the first supercompact cardinal, so that the \SCH\ holds
when I need it.  But let me begin the argument at hand. It suffices to show
that all regular cardinals below $\k$ are preserved. Suppose towards a
contradiction that $\l$ is regular in $V$, but that $\cof(\l)=\g<\l$ in
$V[G]$. The cardinal $\g$ must be regular. Factor $\P$ as
$\P_\g*\Q_\g*\Ptail$, where $\Q_\g$ is nontrivial only if $\g$ is an
allowed stage. The forcing $\Ptail$ is ${\lte}\g$-closed, and so it could
not collapse the cofinality of $\l$ to $\g$. Similarly, the forcing
$\Q_\g$ would not be allowed if it collapsed the cofinality of $\l$.
So I need only check that $\P_\g$ does not collapse the cofinality of
$\l$. Let $\b$ be the supremum of the allowed stages before $\g$; I
really need to check only that $\P_\b$ does not collapse the cofinality
of $\l$. By stripping off the successor stages one by one, since these
cannot collapse cofinalities, I may assume
that $\b$ is a limit of allowed stages. Thus, in particular, $\b$ is a strong
limit cardinal. Since the next stage in a
partial Laver preparation does not occur until beyond the size of the
previous stage forcing, it follows that $\card{\P_\a}<\b$ for all
$\a<\b$. Now, there are two possibilities.  If, on the one hand, $\b$
is regular then, being a regular limit of inaccessibles, it follows
that $\b$ is itself inaccessible. Thus, $\card{\P_\b}=\b$ since I took
a direct limit at stage $\b$, and so $\P_\b$, having size $\b$, is
$\b\plus$-c.c., and consequently unable to collapse the cofinality of
$\l$. If, on the other hand, $\b$ is singular, then $\b<\g$ and also,
since I took an inverse limit at stage $\b$,
$\card{\P_\b}\le\b^{\cof(\b)}=\b\plus$, by the \SCH. Thus, $\P_\b$ is
$\b\plusplus$-c.c., and this is good enough since
$\b<\g<\l\implies\b\plusplus\le\l$, so the cofinality of $\l$ could not
have been collapsed by $\P_\b$. In any case, therefore, I obtain a
contradiction. So neither cardinals nor cofinalities are collapsed
between $V$ and $V[G]$. And this completes the proof.\qed

\theorem Exact Preservation Theorem. Assume that $\k$ is supercompact
in $V$. Then there is a forcing extension $V[G]$ in which $\k$ remains
supercompact, the \GCH\ holds, and over which the supercompactness of
$\k$ is preserved by exactly those $\ltk$-directed closed posets which
preserve the \GCH.

\proof For this proof, I will not actually need to use the Gap Forcing
Theorem \GapForcingTheorem.  I may assume, by forcing if necessary,
that the \GCH\ holds in $V$.  Suppose that $G$ is $V$-generic for the
partial Laver preparation of $\k$ in which stage $\g$ is allowed
provided that it preserves the \GCH\ (cardinals may be collapsed). By
the Partial Laver Preparation Lemma, $\k$ remains supercompact in
$V[G]$.

First, I will argue that the \GCH\ still holds in $V[G]$. Certainly it
still holds at $\k$ and above, since $\P$ has size $\k$. Suppose that
the \GCH\ holds up to $\g<\k$. I can factor $\P$ as
$\P_\g*\Q_\g*\Ptail$, where $\Q_\g$ is trivial unless $\g$ is allowed.
The tail forcing $\Ptail$ adds no subsets to $\g$, and so cannot affect
the \GCH\ at $\g$. The stage $\g$ forcing $\Q_\g$ is only allowed
provided that it preserves the \GCH. So consider $\P_\g$. Let $\b$ be
the supremum of the allowed stages below $\g$. It suffices to show that
the $\P_\b$ preserves the \GCH\ at $\g$.  Without loss of generality,
by stripping off the successor stages one by one, which cannot affect
the \GCH, I may assume that $\b$ is a limit of allowed stages. Now, if
$\b<\g$, then $\P_\b$, being the limit of smaller posets, has
size $\b\plus\leq\g$, and so it cannot destroy the \GCH\ at $\g$. Otherwise,
assume $\b=\g$, so $\g$ is a limit of allowed stages, and therefore a
strong limit cardinal. If $\g$ is regular, then it must be
inaccessible, and so $\P_\g$, using the direct limit, has size $\g$,
and consequently cannot destroy the \GCH\ at $\g$. So assume $\g$ is
singular. By Silver's theorem \cite[Sil74], it suffices to consider the
case that $\cof(\g)=\w$, since otherwise the \GCH, holding below $\g$,
automatically holds at $\g$. But in this case, $2^\g=\g^\w$, and since
the entire forcing $\P$ is countably closed, it cannot affect $\g^\w$.
Thus, it cannot destroy the \GCH\ at $\g$. So $V[G]\sat\GCH$.

Let me now argue that $\k$ is indestructible by any $\ltk$-directed closed
forcing which preserves the \GCH. Suppose $\Q$ is such forcing, and
that $H\of\Q$ is $V[G]$-generic. Fix $\l$ and pick $\theta$ much larger
than both $\l$ and $\card{\Q}$, and a $\theta$-supercompact embedding
$j:V\to M$ such that $j(\ell)(\k)=\Qdot$ and
$\dom(j(\ell))\intersect(\k,\theta]=\emptyset$. Since $H$ preserves the
\GCH\ over $V[G]$, it follows that $H$ also preserves the \GCH\ over
$M[G]$, since $M[G]$ and $V[G]$ agree up to $\theta$, which is much
larger than $H$. The stage $\k$ forcing, therefore, is allowed in
$j(\P)$, and I may continue the usual argument to lift $j$ to
$j:V[G][H]\to M[j(G)][j(H)]$, and then use $j\image\l$ as a seed to
conclude that $\k$ is still $\l$-supercompact in $V[G][H]$.

Finally, it is easy to see that the supercompactness of $\k$ is
destroyed by any $\ltk$-closed forcing which does not preserve the
\GCH, and this is simply because if the \GCH\ holds up to a
supercompact cardinal $\k$, then it must hold everywhere.\qed

Similar arguments establish the next two theorems, whose proofs
I omit, except to say that in the first, one uses a partial Laver
preparation in which $\g$ is allowed when $\Q_\g$ does not collapse
$\g\plus$, and in the second, $\g$ is allowed when $\g$ is closed
under some fixed high jump function $h:\k\to\k$ and $\Q_\g$ {\it does}
collapse $\g\plus$.

\theorem Exact Preservation Theorem. Assume that $\k$ is supercompact
in $V$. Then there is a forcing extension $V[G]$ in which $\k$ remains
supercompact and over which the supercompactness of $\k$ is preserved by
exactly those $\ltk$-directed closed posets which do not collapse
$\k\plus$. Indeed, in $V[G]$ a $\ltk$-directed closed poset, if it
preserves $\k\plus$, will preserve the supercompactness of $\k$; if
it collapses $\k\plus$, it will destroy the measurability of $\k$.

\theorem Exact Preservation Theorem. Assume that $\k$ is supercompact
in $V$. Then there is a forcing extension $V[G]$ in which $\k$ remains
supercompact and over which the supercompactness of $\k$ is preserved by
exactly those $\ltk$-directed closed posets which collapse $\k\plus$.

There is no end to these kinds of theorems. In particular, it is an
easy matter to change the $\k\plus$ in the previous two theorems to
$\k\plusplus$ or $\k\plusplusplus$ and so on: I can get a model
where the supercompactness of $\k$ is preserved by exactly those
$\ltk$-directed closed posets which collapse $\k\plusplusplusplus$ to
$\k\plusplus$, to name but one way of modifying the theorem.  Really it
is only due to the Gap Forcing Theorem \GapForcingTheorem\
that we have some knowledge
about the supercompact embeddings which live in a gap forcing extension
of $V$. In the next section I will introduce a new topic in order
to prove, later, even more powerful exact preservation theorems than these.

\remark Remark on Closure. I have made some mention of forcing notions
which are, variously, $\ltk$-directed closed, $\ltk$-closed, and
$\ltk$-strategically closed. But how much closure do I need for the arguments?
The answer is that I need the forcing $\Q$ to be $\ltk$-directed closed when
I want to argue
that $\Q$ preserves supercompactness; directed closure is used to find
a master condition for $j(\Q)$. I only need $\Q$, however, to be
$\ltk$-closed, or, even weaker, $\ltk$-strategically closed, when I want
to argue that $\Q$ destroys supercompactness, since this is all that the
Gap Forcing Theorem \GapForcingTheorem\
requires. I have stated all the Exact Preservation Theorems
only for $\ltk$-directed closed forcing, for simplicity, but the proofs
go through without modification for strategically closed posets on the
destruction side, provided one allows strategically closed forcing in the
partial Laver preparation itself.

\section Representability

In order to improve the Exact Preservation Theorems, I will now
generalize the beautiful fact of folklore---happily discovered I am
sure by many young set theorists, perhaps like myself, while sipping
coffee in cafes---that for every ordinal $\a<\k\plus$ there is a
function $f:\k\to\k$ which represents $\a$ with respect to every normal
measure on $\k$. Specifically, $j(f)(\k)=\a$ for any such ultrapower embedding.

More generally, now, supposing first that $a\in H(\theta\plus)$, I will say
that $a$ is {\df represented} by the function $f:\k\to V_\k$ with respect to
$\theta$-supercompact embeddings exactly when $$j(f)(\k)=a$$ for any
such embedding. Generalizing still further, to allow for the possibility that
$a$ is much larger than $\theta$, or that
$a$ is a proper class, I officially require only that
$$j(f)(\k)\intersect H(\theta\plus)=a\intersect H(\theta\plus)$$ for
all the $\theta$-supercompact embeddings $j$. This agrees with the first
definition when $a\in H(\theta\plus)$. Thus I require
that $j(f)(\k)$ agree with $a$ as well as the embeddings can be expected to
make it so. I am referring here only to embeddings with critical
point $\k$. Define that $a$ is {\df representable} when there is a
function $f\from \k\to V_\k$ which represents $a$ with respect to
$\theta$-supercompact
embeddings for every $\theta\geq\k$.  The set $a$ is {\df eventually}
representable if there is a function which represents $a$ with respect
to all $\theta$-supercompact embeddings for sufficiently large
$\theta$, and $a$ is {\df frequently} representable if there is a
function which represents $a$ with respect to all $\theta$-supercompact
embeddings, for arbitrarily large $\theta$.

\lemma Folklore Fact. Every ordinal below $\k\plus$ is representable.

\proof This is the happy fact I referred to above.  Certainly any
ordinal below $\k$ is representable; one simply uses a constant
function. Now suppose $\k\leq\a<\k\plus$.  Therefore,
$\a=\ot\<\k,\triangleleft >$ for some relation $\triangleleft$ on
$\k$.  Define $f(\g)=\ot\<\g,\triangleleft\restrict\g>$. Thus, for any
$j:V\to M$ with critical point $\k$, regardless of where $j$ is
defined, it follows that
$j(f)(\k)=\ot\<\k,j(\triangleleft)\restrict\k>= \ot\<\k,\triangleleft
>=\a$, as desired.\qed Fact

But this fact is only the start. What I aim to show now is that there are many
more representable ordinals; indeed, the next two closure
theorems show that the set of representable sets forms a small
set-theoretic universe.

\theorem Ordinal Representation Theorem. The set of (eventually)
representable ordinals is closed under ordinal and cardinal
arithmetic.  Specifically, if $\a$ and $\b$ are (eventually)
representable, then so are the ordinals $\a+\b$, $\a\b$, $\a^\b$, and
the cardinals $\card{\a}$, $\aleph_\a$, $\beth_\a$, $\a\plus$, and
$\card{\a}^{\card{\b}}$. Furthermore, the set of (eventually)
representable ordinals is closed under $\ltek$-suprema.

\proof Suppose that $\a$ and $\b$ are represented by the functions
$f_\a$ and $f_\b$. It is easy to see that $\a+\b$, etc. are represented
by the following functions:  $$f_{\a+\b}(\g)=f_\a(\g)+f_\b(\g)\qquad
      f_{\a\b}(\g)=f_\a(\g)\cdot f_\b(\g)\qquad
      f_{\a^\b}(\g)=f_\a(\g)^{f_\b(\g)}$$
This follows by the absoluteness of ordinal arithmetic between $M$ and
$V$.  For example, if $j$ is $\theta$-supercompact and $\a,\b\leq\theta$, then
$j(f_{\a+\b})(\k)=j(f_\a)(\k)+j(f_\b)(\k)=\a+\b$, as
desired.  The cardinals $\card{\a}$, etc. are represented by the
functions:
$$\vbox{\halign{\hfil#&#\hfil&\qquad\hfil#&#\hfil&\qquad\hfil#&#\hfil\cr
 $f_{\card{\a}}(\g)={}$&$\card{f_\a(\g)}$&
       $f_{\aleph_\a}(\g)={}$&$\aleph_{f_\a(\g)}$&
       $f_{\beth_\a}(\g)={}$&$\beth_{f_\a(\g)}$\cr
       $f_{\a\plus}(\g)={}$&$f_\a(\g)\plus$&
  $f_{\card{\a}^{\card{\b}}}(\g)={}$&$\card{f_\a(\g)}^{\card{f_\b(\g)}}$&\cr}}$$
Again, I need to appeal to absoluteness for these notions between $M$
and $V$. For illustration, consider $\aleph_\a$. Suppose that $j:V\to
M$ is a $\theta$-supercompact embedding. The argument is a little
easier in the case that $\theta$ is very large, for then
$j(f_\a)(\k)=\a$ and consequently
$j(f_{\aleph_\a})(\k)=\aleph_\a^M=\aleph_\a$, as required.  If $\theta$
is less than $\aleph_\a$, then I only need to show that
$j(f_{\aleph_\a})(\k)$ is at least $\theta\plus$. In the case that
$\theta$ is between $\a$ and $\aleph_\a$ then $j(f_\a)(\k)=\a$ and
consequently $j(f_{\aleph_\a})(\k)=\aleph_\a^M$, which is at least
$\theta\plus$. If $\theta$ is less than $\a$, then $j(f_\a)(\k)$ is at
least $\theta\plus$, so $j(f_{\aleph_\a})(\k)=\aleph_{j(f_\a)(\k)}^M$,
which is at least $\theta\plus$. So in any case $\aleph_\a$ and
$j(f_{\aleph_\a})(\k)$ agree up to $\theta\plus$, as is required. The
other cases are similar.

Finally, let me show that the set of (eventually) representable
ordinals is closed under $\ltek$-sups. Suppose that
$\l=\sup_{\a<\k}\l_\a$, where $\l_\a$ is represented by $f_{\l_\a}$.
Let $f_\l(\g)=\sup_{\a<\g}f_{\l_\a}(\g)$, and observe that
$j(f_\l)(\k)=\sup_{\a<\k}j(f_{\l_\a})(\k)$. This agrees with $\l$ up to
$\theta\plus$ since $j(f_{\l_\a})(\k)$ agrees with $\l_\a$ up to
$\theta\plus$.\qed

\theorem Representation Theorem. Every element of $H(\k\plus)$ is
representable. The set of (eventually) representable sets is closed
under elementary set operations. Specifically, if $a$ and $b$ are
(eventually) representable sets, then so are $\{a,b\}$,
$a\minus b$, and $P(a)$; as are the set $\set{x\in a\st \phi(x,b)}$ and
the ordinal $\mu\a[\phi(H(\a\plus),\a,a)]$, when $\phi$ is any
$\Delta_0$ formula. If $a$ is eventually representable, then so is $\union a$.
If $\a$ is an (eventually) representable ordinal,
then $V_\a$ is (eventually) representable. Also, the set of
(eventually) representable sets is closed under $\k$-sequences.

\proof Suppose that $a\in H(\k\plus)$. Thus, $a$ can be coded by some
set $A\of\k$. Let $f_a(\g)$ be the set in $H(\g\plus)$ which is coded
in the same way by $A\intersect \g$, so that $j(f)(\k)$ is the set
coded by $j(A)\intersect \k=A$. That is, $j(f)(\k)=a$, as required.

The closure claims are similar to the previous theorem. One simply uses
the obvious function in each case, and then appeals to absoluteness
between $M$ and $V$. Let me illustrate with $P(a)$. Suppose that $a$ is
represented by the function $f_a$. Define $f_{P(a)}(\g)=P(f_a(\g))$,
and suppose that $j:V\to M$ is a $\theta$-supercompact embedding. Since
$j(f_a)(\k)$ agrees with $a$ on $H(\theta\plus)$, and $M$ is closed
under $\theta$-sequences, it follows that $P(j(f_a)(\k))^M$ agrees with
$P(a)$ on $H(\theta\plus)$, as is required. The other cases are
similar.  To show the closure under $\k$-sequences, suppose that $a_\a$
is represented by $f_\a$ for $\a<\k$. It is easy to see that the
function $f(\g)=\seq<f_\a(\g)\st\a<\g>$ represents the sequence
$\seq<a_\a\st\a<\k>$.\qed

It seems possible that $a$ could be representable but not $\union a$,
because for some $\theta$ perhaps $a$ and $j(f)(\k)$ disagree about
a set which is not in $H(\theta\plus)$ but which has elements in
$H(\theta\plus)$. This is why I only make the claim of eventual
representability for $\union a$ in the theorem.

\lemma Going-Up Lemma. If $f$ represents $a$ with respect to all
$\l$-supercompact embeddings, and $a\in H(\l\plus)$, then $f$
represents $a$ with respect to all $\theta$-supercompact embeddings for
any $\theta>\l$.

\proof Suppose that $j:V\to M$ is a $\theta$-supercompact embedding,
for $\theta>\l$. Let $\mu$ be the $\l$-supercompactness
measure germinated by the seed
$j\image\l$. That is, $X\in\mu\iff j\image\l\in j(X)$. If $\jmu:V\to
M_\mu$ is the ultrapower by $\mu$, it follows that $a\in M_\mu$, and
also that $\jmu=\pi\compose j$, where $\pi$ is the collapse of the seed
hull $X=\set{j(g)(j\image\l)\st g\in V}\elesub M$ (consult \cite[Ham97b]
for elaboration on this seed hull factor method), as illustrated in the
following diagram.
$$\trianglediagram{V}{\jmu}{j}{M_\mu}{\pi^{-1}}{M}$$
Since $\a=\ot(j\image\l\restrict j(\a))$ for all $\a\leq\l$ it follows that
$\l\of X$ and $\l\in X$. Thus, $\pi(a)=a$ and
$\pi(\k)=\k$, and consequently, since $f$ represents $a$ with respect
to $\jmu$, it follows that $a=\jmu(f)(\k)=\pi(j(f))(\k)=\pi(j(f)(\k))$.
Thus $j(f)(\k)=a$ as desired.\qed

\corollary. A set is frequently representable iff it is eventually
representable.

\proof Immediate from the Going-Up Lemma.\qed

\lemma Going-Down Lemma. If $f$ represents $a$ with respect to a
$\theta$-supercompact embedding, then for every $\l<\theta$ there is a
$\l$-supercompact embedding with respect to which $f$ represents $a$.
\ref\GoingDownLemma

\proof Suppose $f$ represents $a$ with respect to $j:V\to M$, a
$\theta$-supercompact embedding. Let $\mu$
be the measure on $P_\k\l$ germinated by $j\image\l$ via $j$, and let
$\jmu:V\to M_\mu$ be the corresponding ultrapower embedding. If $x\in
H(\l\plus)$, it follows as in the previous argument that $\pi(x)=x$,
and consequently $$\eqalign{x\in a\iff&\,x\in j(f)(\k)\cr
         \iff&\,\pi(x)\in\pi(j(f)(\k))\cr
         \iff&\,x\in\jmu(f)(\k).\cr}$$
So $f$ represents $a$ with respect to $\jmu$, a $\l$-supercompact
embedding.\qed

\theorem Enduring Representability Theorem. Suppose that a class $A$ is
represented by the function $f:\k\to V_\k$ in $V$, and that $V[G]$ is a
forcing extension which admits a gap below $\k$. If either (1) $\k$ is
supercompact in $V[G]$, or (2) $A$ is a set and $\k$ is
$2^{\card{\TC(A)}}$-supercompact in $V[G]$, or (3) the \GCH\ holds in
$V$, then $f$ continues to represent $A$ in $V[G]$.
\ref\EnduringRepresentabilityTheorem

\proof Suppose that $j:V[G]\to M[j(G)]$ is a $\theta$-supercompact
embedding in $V[G]$. I need to show, under the various hypotheses, that
$j(f)(\k)$ agrees with $A$ on $H(\theta\plus)$. By the Going-Up Lemma,
if $A$ is a set I may assume that $\theta\leq\card{\TC(A)}$. The
argument is complicated somewhat by the possibility that
${\theta\plus}^V$ may be collapsed by $G$. Suppose that $x\in
H({\theta\plus}^{V[G]})$. I aim to show that $x\in A\iff x\in
j(f)(\k)$. Since $j(f)(\k)\in M\of V$, it suffices to consider only the
case when $x\in V$. It follows that there is some $\l$ such that
$\theta\leq\l<{\theta\plus}^{V[G]}$ and $x\in H(\l\plus)^V$. Since
$\card{\l}=\theta$ in $V[G]$, I may view $j$ as a $\l$-supercompact
embedding. Let $\mu$ be the measure germinated by $j\image\l$ via
$j\restrict V$. By the Local Improvements \LocalImprovements\ of the
Gap Forcing Corollary \GapForcingCorollary, under any of the hypotheses
in 1, 2, or 3, this measure, and therefore also the corresponding
ultrapower map $\jmu:V\to M_\mu$, is in $V$.  Consequently, since $f$
represents $A$ in $V$, it follows that $\jmu(f)(\k)$ agrees with $A$ on
$H(\l\plus)^V$. Since $\mu$ is germinated from $j\image\l$ via $j$, it
follows that $\jmu=\pi\compose j$ where $\pi$ is the collapse of the
seed hull $X=\set{j(g)(j\image\l)\st g\in V}\elesub M$. Since as in the
Going-Up Lemma $\l\of X$ and $\l\in X$, and also
$H(\l\plus)^{M_\mu}=H(\l\plus)^V$, it follows that $\pi(x)=x$ and
$\pi(\k)=\k$.  Now simply compute:
$$\eqalign{x\in A\iff&\;x\in \jmu(f)(\k)\cr
         \iff&\;\pi(x)\in\pi(\jmu(f)(\k))\cr
         \iff&\;x\in\pi(\jmu(f))(\pi(\k))\cr
                 \iff&\;x\in j(f)(\k).\cr}$$
Thus $A$ and $j(f)(\k)$ agree on $x$, as is required.\qed

While the Representation Theorems show that the class of representable
sets forms a small set-theoretic universe, we must keep in mind that
a single function
$f:\k\to V_\k$ represents at most one set, and so the number
of representable sets is at most $2^\k$. Nevertheless, I can make any set
representable by simply collapsing it to $\k$:

\theorem Forcing Representability. Any set can be made representable by
forcing.

\proof Fix any set $a$. First make $\k$ indestructible. Then, collapse
cardinals to $\k$ so that $a\in H(\k\plus)$.
By the Representation Theorem, this makes $a$ representable.\qed

Because cardinals are collapsed, this proof may be unsatisfying.
One easy improvement is to realize that if $\a<\aleph_\a$, I
can make $\aleph_\a$ representable by collapsing only $\a$ to $\k$.
Iterating this, I can make $\aleph_{\aleph_\a}$ representable by
collapsing only $\a$.  But still cardinals are collapsed.
The next theorem shows how to avoid this, and add {\it any}
set, while collapsing no cardinals above $\k$, to the collection of
{\it eventually} representable sets.

\theorem Forcing Eventual Representability. If $\k$ is supercompact,
then any set can be made eventually representable by forcing which
preserves the supercompactness of $\k$, does not collapse cardinals
above $\k$, and preserves all previously representable and eventually
representable sets.
\ref\ForcingEventualRepresentability

\proof First, using Exact Preservation Theorem \EPCard,
I may assume that the supercompactness of $\k$ is
indestructible by $\ltk$-directed closed
forcing notions which collapse neither cardinals nor
cofinalities. Since the forcing to accomplish this was a partial Laver
preparation, it preserves the supercompactness of $\k$ and admits a gap
below $\k$. Thus, by the Enduring Representability Theorem
\EnduringRepresentabilityTheorem, it
preserves representability. Also, it collapses no cardinals above
$\k$.  Now I will make the set $a$ easily definable by coding it into
the continuum function above $\k$. The usual way of doing this,
however, collapses cardinals in the case that the \GCH\ fails, but I
need not worry, since by \cite[Sol74] the \SCH\ holds above any
supercompact cardinal, and this will be enough for my argument. By the
\SCH, if $\l$ is a singular strong limit above
$\k$ it follows that $2^\l=\l\plus$. I may, therefore, add subsets to
$\l\plus$ without collapsing cardinals. I may assume that $a$ is a set of
ordinals below some $\d$ (to make decoding easier, I may even assume by further
coding that $a$
consists entirely of successor ordinals, except for it's maximum element).
Let $\<\l_\a\st \a<\d>$ enumerate the first
$\d$ many singular strong limits above $\k$. I will design a forcing
notion which will ensure that
$2^{\l_\a\plus}=\aleph_\b$ for $\b$ even or odd, respectively, according
to whether $\a\in a$ or not. Let $\P$ be the reverse Easton support iteration
which at stage $\l_\a$ forces to add $\aleph_{\b+1}$ many subsets to $\l\plus$,
if $\b$ was even and I want it to be odd or vice versa. By a $\Delta$-system
argument (and this is where I use that $2^{\l_\a}=\l_\a\plus$) the stage
$\l_\a$ forcing $\Q_{\l_\a}$ is $\l_\a\plusplus$-c.c., and, since it is
also ${\lte}\l_\a$-closed, each stage of this iteration preserves all cardinals
and cofinalities. Next, a factor argument like that in Theorem \EPCard\
establishes, since the \SCH\ holds above $\k$,
that the entire iteration preserves all cardinals and cofinalites.
This forcing makes the set $a$ concretely definable from the
continuum function. Let $f(\g)$ be the set that is obtained by
running the decoding of this information for the singular strong limits
above $\g$. Then, if $j:V\to M$ is a $\theta$-supercompact embedding
and $\theta$ is large enough so that $M$ has the same continuum
function as $V$ as high as any coding that I performed, then $j(f)(\k)$
will perform the same decoding in $M$ as I did in $V$, and hence
$j(f)(\k)=a$. Thus, $f$ eventually represents the set $a$. This forcing
preserves all previously
(eventually) representable sets because it admits a gap below $\k$
while preserving the supercompactness of $\k$.\qed

\theorem Forcing Frequent Representability Theorem. If $\k$ is
supercompact, then any class of cardinals can be made frequently
representable by
forcing which preserves the supercompactness of $\k$, does not collapse
cardinals above $\k$, and preserves all previously representable and
eventually representable sets and classes.

\proof Just perform the same coding as in the previous theorem. This
time, however, one cannot get above all the coding. Rather, one can
find arbitrarily high $\theta$ which are closure points of the coding
in the sense that $A\intersect\theta=A\intersect H(\theta\plus)$
is coded below $\theta$.
So in this case I conclude only frequent representability.\qed

\theorem Forcing Cardinal Representability Theorem. Assume the \SCH\ holds,
$2^\k=\k\plus$, and $\g^\k=\g$, where $\k$ is supercompact.
Then $\g$ can be made representable without collapsing cardinals.
\ref\ForcingCardinalRepresentabilityTheorem

\proof Using Theorem
\EPCard, I can ensure that the supercompactness of $\k$ is preserved by any
$\ltk$-directed closed forcing which collapses neither cardinals nor
cofinalities. And this can be done while collapsing neither cardinals
nor cofinalities. Now
simply force $2^\k=\g$, making $\g$ representable by the Representation
Theorem. The preparatory
forcing ensures that this will preserve the supercompactness of $\k$,
and neither cardinals nor cofinalities are collapsed. \qed

\corollary. If $2^\k=\k\plus$, where $\k$ is supercompact, and $\g^\k=\g$,
then $\g$ can be made
representable by forcing which does not collapse cardinals above $\k$.

\proof The same argument as the previous theorem. We don't need the
\SCH\ if we don't mind collapsing a few cardinals below $\k$.\qed

\question. Assume the \GCH. Can any set can be made representable by
forcing which does not collapse cardinals?
\ref\Representable

Perhaps this question will be answered by first making the desired
set definable. The standard trick of coding the set into the
continuum function---one makes the \GCH\ hold or fail at
successive cardinals in such a way so as to code the given set---does not collapse
cardinals or cofinalities when the \GCH\ holds; my proof of the
Forcing Eventual Representability Theorem \ForcingEventualRepresentability\
shows that this can be done assuming only the \SCH\ holds. A very interesting
open question is whether there is some clever way of showing just in \ZFC\ that
any set can be made definable without collapsing cardinals or cofinalities. The
difficulty, however, with this whole approach as an attack on Question
\Representable\ is that it is not clear that definability will give
representability, even if there is some simple coding involved, since
in a sense representability is local---we seem to need the information
about $a\intersect H(\theta\plus)$ to be coded into $H(\theta\plus)$,
so that a $\theta$-supercompact embedding has access to it. This is why
the Forcing Representability theorems conclude only the eventual or
frequent representability of $a$. Fortunately for the results of this
paper, this amount of representability goes a long way.

\section Separating the Superdestructibility Hierarchy

While the levels of the supercompactness hierarchy become steadily
stronger as one moves upward, this is not true of the
superdestructibility hierarchy (remember that a cardinal $\k$ is
superdestructible at $\theta$ if any $\ltk$-closed forcing which adds a
subset to $\theta$ destroys the $\theta$-supercompactness of $\k$). The
essential reason for this is that, with a larger $\theta$,
superdestructibility requires a stronger property to be destroyed by a
larger class of forcing notions.  So there is no clear implication
either upwards or downwards. In fact, as I will prove in this section,
I can turn superdestructibility on, then off, and then on again up through
the hierarchy, in almost any conceivable pattern, making, for
example, a supercompact cardinal $\k$ superdestructible at $\k\plus$
and $\k\plusplusplus$ but not at $\k\plusplus$ or
$\k\plusplusplusplus$, and so on.  Let me begin by separating just two
levels of superdestructibility.

\theorem Separation Theorem. Suppose that $\k\leq\l<\theta$, where $\k$
is a supercompact cardinal, $\l$ and $\theta$ are regular, and the
\GCH\ holds.  Then, while collapsing neither cardinals nor cofinalities
above $\k$, one can make $\k$ superdestructible at $\l$ but not at
$\theta$, and, vice versa, superdestructible at $\theta$ but not at
$\l$.

\proof I may assume, by the Forcing Cardinal Representability Theorem
\ForcingCardinalRepresentabilityTheorem, that
$\l$ is represented by some function $f:\k\to\k$. Let me point out that
I will not actually use much of the \GCH. So far I have only used that
$2^\k=\k\plus$ and $\l^\k=\l$, when $\l>\k$, to make $\l$ representable. And
this may even destroy the \GCH\ at $\k$. Later, to build the second model,
I will use that $2^\l\leq\theta$, and that is all of the \GCH\ that I will
assume. In the two
constructions below I will actually prove much more than I claimed.
Suppose that $\ell$ is a Laver function.

\quiet\lemma The First Model. There is a forcing extension over which the
supercompactness of $\k$ is preserved by the $\ltk$-directed closed
posets which do not add a subset to $\l$, and over which the
$\l$-supercompactness of $\k$ is destroyed by those which do.

\proof Suppose that $G$ is $V$-generic for the partial Laver
preparation $\P$ of $\k$ in which stage $\g$ is allowed provided that it does
not add subsets to $f(\g)$.  By the Partial Laver Preparation Theorem,
$\k$ remains supercompact in $V[G]$. The usual argument shows that the
supercompactness of $\k$ is preserved by any $\ltk$-directed closed
forcing which does not add subsets to $\l$.

Let me now show that in $V[G]$, the $\l$-supercompactness of $\k$ is
destroyed by any $\ltk$-closed forcing which adds a subset to $\l$.
That is, I will show that $\k$ is superdestructible at $\l$ in $V[G]$.
If this is not true, then there is some $V[G]$-generic $H\of\Q$,
where $\Q$ is $\ltk$-closed, which adds a new set $B\of\l$, such that
$\k$ remains $\l$-supercompact in $V[G][H]$.  Let $j:V[G][H]\to
M[j(G)][j(H)]$ be a $\l$-supercompact embedding. By the Enduring
Representability Theorem \EnduringRepresentabilityTheorem,
I know $j(f)(\k)=\l$.  Since, by the closure
of the embedding, $B\in M[j(G)][j(H)]$, it follows that $B\in M[j(G)]$
by the closure of the $j(\Q)$ forcing.  Also, the stage $\k$
forcing in $j(G)$ is not allowed to add new subsets to $j(f)(\k)=\l$,
and I know the stages after $\k$ do not begin until after $\l$, and
consequently they are ${\lte}\l$-closed. It follows that $B\in
M[G]$, and so, by the Gap Forcing Theorem \GapForcingTheorem, $B\in V[G]$, contradicting
that $B$ was newly added by $H$.\qed First Model

\quiet\lemma The Second Model. There is a forcing extension in which
$\k$ remains supercompact and the $\l$-supercompactness of $\k$ becomes
fully indestructible by any $\ltk$-directed closed forcing, but over
which any $\ltk$-closed forcing which adds a subset to any $\theta\geq 2^\l$
will destroy the $\theta$-supercompactness of $\k$.

\proof I may assume by the High Jump Theorem \HighJumpTheorem\ that there
is a high jump function $h$ for $\k$.
Let $\P$ be the partial Laver preparation of $\k$ in which stage
$\g$ is allowed provided that, first, it destroys the $f(\g)$-supercompactness
of $\g$, and, second, $\g$ is a closure point of the
functions $f$ and $h$. Suppose that $G$ is $V$-generic for $\P$. I
will show that $V[G]$ has the desired properties. By the Partial Laver
Preparation Lemma \PartialLaverPreparationLemma, I know that $\k$ remains supercompact in $V[G]$.
And furthermore, neither cardinals nor cofinalities above $\k$ are
collapsed, because the forcing $\P$ is $\k$-c.c.

Let me now prove that the $\l$-supercompactness of $\k$ is fully
indestructible in $V[G]$. Suppose that $H\of\Q$ is $V[G]$-generic for the
$\ltk$-directed closed forcing $\Q=\Qdot_G$, and, towards a
contradiction, that $\k$ is not $\l$-supercompact in $V[G][H]$. Fix
$\d\muchgt\l$ and a $\d$-supercompact embedding $j:V\to M$ such that
$j(\ell)(\k)=\Qdot$ and $\dom(j(\ell))\intersect (\k,\d]=\emptyset$.
Notice that the stage $\k$ forcing in $j(\P)$ will be $\Qdot$, provided
that $\k$ is allowed. And, since $G*H$ is $M$-generic for $\P*\Qdot$, I
know that $\k$ will be allowed if it is not $\l$-supercompact in
$M[G][H]$, since by representability $j(f)(\k)=\l$. But since I assumed
that $\k$ was not $\l$-supercompact in $V[G][H]$, it follows that $\k$
is not $\l$-supercompact in $M[G][H]$, since $M$ and $V$ agree up to
$\d$. Thus, the stage $\k$ forcing in $j(\P)$ is allowed, and so
$j(\P)$ factors as $\P*\Qdot*\Ptail$, where $\Ptail$ is
${\lte}\d$-closed, since the high jump function jumps over $\d$.  I can
therefore proceed as usual to lift the embedding to $j:V[G][H]\to
M[j(G)][j(H)]$ by forcing to add $\Gtail$ and $j(H)$. Again, using
$j\image\l$ as a seed, I generate a normal fine measure on $P_\k\l$
which can not have been added by the tail forcing, and consequently
must be in $V[G][H]$, contradicting our assumption that $\k$ was not
$\l$-supercompact there.

Finally, let me prove that any $\ltk$-closed forcing which adds a
subset to $\theta\geq 2^\l=2^{\l^\ltk}$ will destroy the
$\theta$-supercompactness of $\k$. Suppose towards a
contradiction that $\k$ remains $\theta$-supercompact in $V[G][H]$,
where $H\of\Q$ is $V[G]$-generic and adds a subset $B\of\theta$.  Then
there is a $\theta$-supercompact embedding $j:V[G][H]\to
M[j(G)][j(H)]$.  By the Enduring Representability Theorem
\EnduringRepresentabilityTheorem, I know
$j(f)(\k)=\l$. If $\k$ is allowed, then $\k$ is
not $\l$-supercompact in $M[G][g]$, where $g$ is the stage $\k$ generic
of $j(G)$. But the additional forcing to $M[j(G)][j(H)]$ is
${\lte}\theta$-closed, since $j(h)(\k)>\theta$,
and so $\k$ is not $\l$-supercompact in
$M[j(G)][j(H)]$. But $\k$ is $\l$-supercompact in $V[G][H]$, so by
coding a measure on $P_\k\l$ with a subset of $\theta$ and using the closure of
$j$, we see that $\k$ is $\l$-supercompact in $M[j(G)][j(H)]$, a
contradiction. Thus, $\k$ must not be allowed. In this case, I know
that $B\in M[G]$ by the closure of the tail forcing, and so, by the
Gap Forcing Theorem \GapForcingTheorem,
$B\in V[G]$, contradicting our assumption that
$B$ was newly added by $H$. So in any case I get a contradiction.\qed
Second Model

This completes the proof of the theorem.\qed

Let me prove next a great generalization of the previous theorem.

\theorem Superdestruction Separation Theorem. Suppose that $A$ is a
representable class of cardinals each with cofinality at least $\k$,
and that the \GCH\ holds. Then there is a forcing extension $V[G]$,
preserving all cardinals, cofinalities, and the \GCH, in which:
\points 1. If $\l\notin A$, then $\k$ is superdestructible at $\l$; any
           $\ltk$-closed forcing which adds a subset to $\l$ destroys the
           $\l$-supercompactness of $\k$.\cr
    2. If $\l\in A$, then $\k$ is not superdestructible at $\l$.
       Indeed, any $\ltk$-directed closed forcing which adds no
       bounded sets to $\l$ and which preserves all cardinals,
       cofinalities and the \GCH\ will preserve the
       $\l$-supercompactness of $\k$.\cr

\proof Suppose that $A$ is represented by the function $f$. By the High
Jump Theorem \HighJumpTheorem, I may assume that there is a high jump function
$h$ for $\k$. Let $\P$ be
the partial Laver preparation of $\k$ in which stage $\g$ is allowed provided
first, that $\Q_\g$ adds no bounded
subsets to some $\l\in f(\g)$;
second, that $\g$ is not $\l$-supercompact, for this
same $\l$, in $V^{\P_\g*\Q_\g}$; third, that $\Q_\g$ preserves all
cardinals, cofinalities, and the \GCH; and finally, fourth, that $\g$
is a closure point both of $f$ and the high jump function $h$, in the
sense that $f\image\g\of V_\g$ and $h\image\g\of\g$. This defines the
forcing $\P$. Now suppose that $G\of\P$ is $V$-generic. Let me prove
that $V[G]$ has the desired properties. I know by the Partial Laver
Preparation Lemma \PartialLaverPreparationLemma\
that $\k$ remains supercompact in $V[G]$.
Furthermore, by the argument of Theorem \EPCard, $\P$ preserves all
cardinals, cofinalities, and the \GCH.

Let me prove that 1 holds. Suppose $\l$ is not in $A$, but that $\k$
remains $\l$-supercompact in $V[G][H]$ where $H\of\Q$ is $V[G]$-generic
for the $\ltk$-closed forcing $\Q$ which adds a subset $B\of\l$. Then
there must be a $\l$-supercompact embedding $j:V[G][H]\to
M[j(G)][j(H)]$ witnessing this. Necessarily, $B\in M[j(G)][j(H)]$.
Because I only forced at stages which were closure points of the
high jump function, it follows that there is no forcing in $j(\P)$ in
the interval $(\k,\l]$. Suppose, momentarily, that
there is no forcing at stage $\k$ in $j(\P)$. In this case,
$j(G)=G*\Gtail$ and by the previous remarks $\Gtail$ is
${\lte}\l$-closed. Consequently, by closure considerations, $B\in M[G]$,
and since $M\of V$ by the Gap Forcing Theorem \GapForcingTheorem\
it follows that $B\in
V[G]$, a contradiction. Thus there must be forcing at
stage $\k$. In this case $j(G)=G*g*\Gtail$, where $g\of\Q_\k$ is the
stage $\k$ forcing, with the corresponding $\l'\in j(f)(\k)$ to which
$g$ adds no bounded sets, and such that $\k$ is not $\l'$-supercompact
in $M[G][g]$. Since every element of $A$ has cofinality at least $\k$, I
may assume $\cof(\l')\geq\k$.
Again I know that $\Gtail$ is ${\lte}\l$-closed. Therefore
$B\in M[G][g]$. Also, since $g$ was allowed I know that it preserved
all cardinals, cofinalities, and the \GCH.  By the closure of the
forcing and of $j$ I know that
$P(\l)^{M[G][g]}=P(\l)^{M[j(G)][j(H)]}=P(\l)^{V[G][H]}$. Thus, by
coding measures from $V[G][H]$ on $P_\k\b$ into subsets of $\l$, I conclude
that $\k$ is
$\b$-supercompact in $M[G][g]$ for every $\b<\l$ with $\cof(\b)\ge\k$.
It follows that $\l\le\l'$. Since $\l\notin A$ and $\l'\in j(f)(\k)$ it
follows that $\l\not=\l'$, since $j(f)(\k)$ and $A$ agree up to
$\l\plus$. Consequently $\l<\l'$. Thus, since $g$ adds no bounded sets
to $\l'$, I conclude that $B\in M[G]$, and, as in the first case, that
$B\in V[G]$, a contradiction. So 1 is proved.

To prove 2, suppose that $\l\in A$ and that $H\of\Q$ is $V[G]$-generic,
where $\Q$ is $\ltk$-closed, adds no bounded sets to $\l$, and
preserves all cardinals, cofinalities, and the \GCH. Towards a
contradiction, suppose that $\k$ is not $\l$-supercompact in
$V[G][H]$.  Fix $\theta\muchgt\l$ and, in $V$, a $\theta$-supercompact
embedding $j:V\to M$, such that $j(\ell)(\k)=\dot\Q$ and
$\dom(j(\ell))\intersect(\k,\theta]=\emptyset$.  Since $\k$ is not
$\l$-supercompact in $V[G][H]$, it is also not supercompact in
$M[G][H]$, and so $\k$ is an allowed stage of $j(\P)$.
I may therefore employ the usual lifting
argument to lift to $j:V[G][H]\to M[j(G)][j(H)]$ in a forcing extension
$V[G][H][\Gtail][j(H)]$. As usual, using $j\image\l$ as a seed, I obtain a
measure $\mu$ on $P_\k\l$ which must be in $V[G][H]$, contrary to our
assumption.\qed

\section Exact Preservation As You Like It

In the previous sections we discovered a few major
landmarks in the unexplored region between indestructibility
and superdestructibility. In this section I will point out a
great mountain range, spanning the continent. Specifically, after
proving the Exact Preservation Theorems, more
powerful than the ones in \S 3, my arguments will culminate in the `As
You Like It' Theorem, the title theorem of this paper.
Recall that a set of ordinals is {\df fresh} over $V$ when every proper initial
segment of it is in $V$, but the set itself is not in $V$. Let us begin.

\theorem Exact Preservation. Suppose that $\k$ is supercompact in $V$
and that $A$ is a class of cardinals. Then there is a forcing extension
$V[G]$, obtained without collapsing cardinals above $\k$,
in which $\k$ remains supercompact and over which the
supercompactness of $\k$ is preserved by exactly those $\ltk$-directed
closed posets which collapse no cardinal in $A$.

\proof Again by forcing if necessary I may assume that $A$ is
frequently represented by some function $f$, and that there is a high jump
function $h$. Now suppose $G$ is
$V$-generic for the partial Laver preparation of $\k$ in which $\g$ is allowed
if $\Q_\g$ collapses no cardinal in $f(\g)$ and $\g$ is closed under
the high jump function. The usual arguments establish that the
supercompactness of $\k$ in $V[G]$ is preserved by any forcing which
collapses no cardinals in $A$. Let me now prove it is destroyed by
those which do. Suppose $H\of\Q$ is $V[G]$-generic for $\ltk$-closed
forcing $\Q$ which collapses a cardinal $\l\in A$, but that $\k$
remains $\theta$-supercompact, for some large $\theta\geq\l$ such that
$f$ represents $A$ with respect to the witness
embedding $j:V[G][H]\to M[j(G)][j(H)]$. (By the Going-Down Lemma
\GoingDownLemma\ I may in fact take $\theta=\l$ here.)
The former cardinal $\l$ must
be collapsed in $M[j(G)][j(H)]$, and since the high jump function jumps
over $\theta$, it must be that $\l$ is collapsed at stage $\k$ in $j(G)$.
But $j(f)(\k)$ and $A$ agree on $\l$, and so the stage $\k$ forcing
would not be allowed if it collapsed $\l$, a contradiction.\qed

\theorem Exact Preservation. Suppose that $\k$ is supercompact in $V$
and $A$ is a class of cardinals. Then there is a forcing extension
$V[G]$, obtained without collapsing cardinals above $\k$,
in which $\k$ remains supercompact and over which the
supercompactness of $\k$ is preserved by exactly those $\ltk$-directed
closed posets which collapse a cardinal of $A$.

\proof Again, by forcing if necessary, I may assume that $A$ is frequently
represented by $f$, and that there is a high jump function $h$. For this
argument, let stage $\g$ be allowed
when $\Q_\g$ collapses a cardinal in $f(\g)$ and $\g$ is closed
under $h$. The preservation arguments go through as
usual. Suppose, conversely, that $H\of\Q$ is $V[G]$-generic for
$\ltk$-closed forcing which collapses no cardinal of $A$, and
that $\k$ remains $\l$-supercompact for some large $\l\geq\card{\Q}$
such that
$f$ represents $A$ with respect to a $\l$-supercompactness  embedding
$j:V[G][H]\to M[j(G)][j(H)]$. Since the high jump function jumps over
$\l$, it follows that $H\in M[G][g]$, where $g$ is the stage $\k$
forcing (the forcing $g$ is nontrivial since $H\notin M[G]$).
Observe that
$V[G][g]=V[G][H]$, and so the $g$ forcing is equivalent over $V[G]$ to
a forcing of size $\leq\l$.  This isomorphism must lie in $M[G]$, and
so $g$ is unable to collapse cardinals above $\l$ in $M[G]$. But $g$
must collapse a cardinal in $j(f)(\k)$, since it was allowed. Thus,
over $M[G]$, and hence also over $V[G]$, the generic $g$
collapses a cardinal in $A$ below $\l$.
This contradicts that $H$ collapsed no such cardinals over $V[G]$.\qed

\theorem Exact Preservation. Suppose that $\k$ is supercompact and $A$
is a class of regular cardinals. Then there is a forcing extension
$V[G]$, obtained without collapsing cardinals above $\k$,
over which the supercompactness of $\k$ is preserved by exactly
those $\ltk$-directed closed posets which add a fresh subset, over
$V[G]$, to a cardinal in $A$.

\proof First, by the High Jump Theorem \HighJumpTheorem, I may assume that
there is a high jump function $h$ for $\k$. Also, by forcing
if necessary, I may assume that $A$ is frequently representable, by some
function $f$. Suppose now that $G$ is generic for the partial Laver
preparation $\P$ of $\k$ in which stage $\g$ forcing is allowed if first,
it adds a fresh subset over $V^{\P_\g}$ to some $\l\in f(\g)$,
and second, $\g$ is closed under $h$. I know by the Partial Laver Preparation
Lemma \PartialLaverPreparationLemma\ that $\k$ remains
supercompact in $V[G]$ and the usual preservation argument shows that
the supercompactness of $\k$ is preserved by any $\ltk$-directed closed
poset which adds a fresh subset to an element of $A$.

Now suppose $H\of\Q$ is $V[G]$-generic for $\ltk$-closed $\Q$ and $H$
does not add a fresh set, over $V[G]$, to any element of $A$. I would
like to show $\k$ is no longer supercompact. Suppose, towards a
contradiction, that $j:V[G][H]\to M[j(G)][j(H)]$ is a $\l$-supercompact
embedding, where $\card{\Q}\leq\l$, and $\l$ is such that $f$
represents $A$ with respect to $j$. Thus, $H\in M[j(G)][j(H)]$, and
since the high jump function $h$ jumps over $\l$, I know therefore that
$H\in M[G][g]$, where $g$ is the (perhaps trivial) stage $\k$
forcing.  Since $H\notin V[G]$, it must be also that $H\notin M[G]$ and so
$g$ is actually nontrivial.
Notice that $V[G][H]=V[G][g]$, and consequently, below a condition, the
$g$ forcing is isomorphic, in $V[G]$, to the $H$ forcing below a condition.
Since this isomorphism has size less than or equal to $\l$, it
must lie in $M[G]$, and consequently, by the chain
condition, $g$ cannot add a fresh subset to any regular cardinal in
$M[G]$ above $\l$. But since $g$ was allowed, it must have added, over
$M[G]$, a fresh subset $B\of\z$ for some $\z\in j(f)(\k)$. By the
previous observation it follows that $\z\leq\l$, and consequently
$\z\in A$. Since $B\notin M[G]$ and $\z\leq\l$ it follows that $B\notin
V[G]$, and so $H$ has added, over $V[G]$, a fresh subset to an element
of $A$, contradicting my assumption.\qed

\theorem Exact Preservation. Suppose that $\k$ is supercompact in $V$
and that $A$ is a frequently representable
class of cardinals at which the
\GCH\ holds.  Then there is a forcing extension $V[G]$, obtained without
collapsing cardinals or disturbing the \GCH\ above $\k$, in which $\k$
remains supercompact, and over which the supercompactness of $\k$ is
preserved by exactly those $\ltk$-directed closed posets which preserve
the \GCH\ at the cardinals of $A$.

\proof Again assume $A$ is frequently represented by $f$, and that there is
a high jump function $h$. This time
let the stage $\g$ forcing be allowed when it preserves the \GCH\ at
every cardinal
in $f(\g)$ and $\g$ is closed under the high jump function. The usual
arguments establish that the supercompactness of $\k$ is preserved over
$V[G]$ by $\ltk$-directed closed posets which preserve the \GCH\ at the
cardinals of $A$. Suppose, conversely, that $H\of\Q$ violates the
\GCH\ at $\l\in A$, but that $\k$ remains $\theta$-supercompact with
embedding $j:V[G][H]\to M[j(G)][j(H)]$, where $\theta\geq\l$ is
such that $f$ represents $A$ with respect to the $\theta$-supercompact
embeddings. Because the high jump function
jumps over $\l$, only the stage $\k$ forcing could ruin the \GCH\ at
$\l$ in $M[j(G)][j(H)]$, but it is in precisely this case that it is not
allowed, a contradiction.\qed

\theorem Exact Preservation. Suppose that $\k$ is supercompact in $V$ and
that $A$ is a frequently representable class of cardinals at
which the
\GCH\ holds. Then there is a forcing extension $V[G]$, obtained without
collapsing cardinals or disturbing the \GCH\ above $\k$, in which $\k$
remains supercompact and over which the supercompactness of $\k$ is
preserved by exactly those $\ltk$-directed closed posets which destroy
the \GCH\ at a cardinal of $A$.

\proof Assume $A$ is frequently represented by $f$, and that $h$ is a high
jump function. Let the stage $\g$ forcing be allowed
when it destroys the \GCH\ at some element of $f(\g)$ and $\g$ is
closed under the high jump function. Again the usual lifting arguments
show that the supercompactness of
$\k$ in $V[G]$ is preserved by any $\ltk$-directed closed forcing which destroys
the \GCH\ at some element of $A$.  Conversely, suppose $H\of\Q$ is
$V[G]$-generic for $\ltk$-closed forcing $H$ which does not destroy the
\GCH\ at any cardinal in $A$, but that $\k$ remains $\l$-supercompact
with embedding $j:V[G][H]\to M[j(G)][j(H)]$ for some large $\l\geq\card{\Q}$
which works with $f$. Again the high jump function jumps
over $\l$, so $H\in M[G][g]$ where $g$ is the stage $\k$ forcing in
$j(\P)$; since $H\notin M[G]$ this forcing is nontrivial.
Observe that $V[G][g]=V[G][H]$, and so the $g$ forcing is equivalent in
$V[G]$ to forcing of size at most $\l$, and the isomorphism must be in
$M[G]$. Thus, over $M[G]$, the generic $g$ does not affect the \GCH\
above $\l$. Since it was allowed, it must have destroyed the \GCH\ at
some element of $j(f)(\k)$ below $\l$, and since $f$ represents $A$,
this element must be in $A$. Thus, over $V[G]$, the generic $g$, and
hence also $H$, destroyed the \GCH\ at an element of $A$, contrary to
my assumption on $H$.\qed

The previous theorems display the power of the Gap Forcing Theorem
\GapForcingTheorem\ to severely limit the sort of supercompactness
embeddings which can exist in a gap forcing extension. All the Exact
Preservation Theorems, however, are special cases of, and follow as
immediate corollaries to, my next theorem,
the `As You Like It' Theorem, which asserts that one can
tailor the universe, by forcing, so that nearly any desired class of posets
will preserve the supercompactness of $\k$, and the others destroy it.
It therefore encompasses all of the particular properties in the Exact
Preservation Theorems, and unifies their proofs.

I will now make two key definitions.
Suppose $\set{\Q\st\phi(\Q,\g,A,G)}$ is a class of
$\ltg$-directed closed posets defined using the formula $\phi$ and, as
parameters, a cardinal $\g$, a class $A$, and a set $G$ such that
$\card{G}\leq\g$. I will say that this class, or the formula $\phi$,
is {\df local} provided that, in any model of set theory,
it can be decided whether a given poset
$\Q\in H(\l\plus)$ is a member of the
class by consulting only $H(\l\plus)$: that is,
first, the truth of $\phi(\Q,\g,A,G)$ does not depend fully on $A$ but rather
only on $A\intersect H(\l\plus)$, and furthermore, second, that the
truth of $\phi(\Q,\g,A,G)$ is absolute to any other model with the same
$H(\l\plus)$. The formula $\phi$, with parameters, {\df respects the
equivalence of forcing} iff in any model of set theory, whenever
$\phi(\Q)$ holds, and $\Q$ and $\Q'$ have isomorphic complete boolean
algebras, then $\phi(\Q')$ also holds; also, $\phi(\Q)$ holds just in
case for densely many $b\in\Q$, $\phi(\Q_b)$ holds, where $\Q_b$
denotes the part of the poset $\Q$ below the condition $b$.  It follows
that $\phi(\Q)$ holds just in case $\phi(\Q_b)$ holds for every
$b\in\Q$. For example, the formulas ``$\Q$ preserves every cardinal in
$A$,'' ``$\Q$ preserves the \GCH\ at the cardinals of $A$,'' and
``$\Q$ adds a fresh subset to an element of $A$'' are all
local definitions which respect the equivalence of forcing: if
$\Q\in H(\l\plus)$, then, since all the relevent names in question,
for the collapsing
functions or the fresh sets, are also in $H(\l\plus)$, it follows that
any other model with the same $H(\l\plus)$ will agree on $\phi(\Q)$.

\theorem The `As You Like It' Theorem. {The class of $\ltk$-directed
closed posets which preserve the supercompactness of $\k$ can be made
by forcing to be defined by any pre-selected local formula which
respects the equivalence of forcing.  \smallskip\noindent More
precisely: suppose that $\k$ is supercompact in $V$ and that
$\phi$ is any local formula you like, with class
parameter $A$, which respects the equivalence of forcing. Then one can
force to a model $V[G]$ where $\k$ remains supercompact, and where, for
any $\ltk$-directed closed poset $\Q$ in $V[G]$:}
\points 1. If $\phi(\Q,\k,A,G)$ holds, then $\Q$ preserves the
       supercompactness of $\k$.\cr
        2. If $\phi(\Q,\k,A,G)$ fails, then, below a condition,
        $\Q$ destroys the supercompactness of $\k$.\cr

\proof Assume that $\k$ is supercompact. I may, by the High Jump Theorem
\HighJumpTheorem, suppose also that there is a high jump function $h$ for $\k$.
Furthermore, by
forcing if necessary, I may assume that $A$ is frequently represented
by some function $f$.  Let $\ell$ be a Laver function. I may assume
that every point in $\dom(\ell)$ is a closure point of the high jump
function $h$, and also of $f$.  The forcing $\P$ will be the partial
Laver preparation of $\k$ in which, at the very first stage,
to avoid triviality,
I add a Cohen real, and then, at subsequent stages $\g$, the forcing
$\Q_\g=\ell(\g)_{G_\g}$ is allowed provided that
$V[G_\g]\sat\phi(\Q_\g,\g,f(\g),G_\g)$. Thus, after the first stage, I
perform the Laver preparation exactly when the Laver function
$\ell$ hands me a poset which satisfies the formula $\phi$ in the
appropriate model.

Suppose now that $G\of\P$ is $V$-generic, and that $\Q$ is a
$\ltk$-directed closed poset in
$V[G]$. I know that $\k$ is supercompact in $V[G]$ by the Partial Laver
Preparation Lemma \PartialLaverPreparationLemma. It remains to prove the other two properties.

\lemma. If $\phi(\Q,\k,A,G)$ holds in $V[G]$, then $\Q$ preserves the
supercompactness of $\k$.

\proof As in the previous theorems, the usual Laver argument adapts to
this circumstance. Suppose
$H\of\Q$ is $V[G]$-generic. Fix $\l\geq\card{\Q}$ and $\theta\muchgt\l$,
and let $j:V\to M$ be a $\theta$-supercompact embedding such that
$j(\ell)(\k)=\Qdot$ and $\dom(j(\ell))\intersect (\k,\theta]=\emptyset$. I
must argue that the stage $\k$ forcing is allowed. Since
$\phi(\Q,\k,A,G)$ holds in $V[G]$, it also holds in $M[G]$, since
$\phi$ is local and $H(\l\plus)^{V[G]}=H(\l\plus)^{M[G]}$. Also, by
representability, $A\intersect H(\theta\plus)=j(f)(\k)\intersect
H(\theta\plus)$, and, since $\phi$ depends not fully on $A$ but only on
$A\intersect H(\theta\plus)$, it follows that $\phi(\Q,\k,j(f)(\k),G)$ holds
in $M[G]$. Thus, the stage $\k$ forcing is allowed. The forcing
$j(\P)$, therefore, factors as $\P*\Q*\Ptail$, where $\Ptail$ is
${\lte}\theta$-closed in $M[G][H]$, and hence also in $V[G][H]$. Thus, as usual,
I can force to add $\Gtail\of\Ptail$ generic over $V[G][H]$, and lift
the embedding to $j:V[G]\to M[j(G)]$ where $j(G)=G*H*\Gtail$. After
this, I can also lift the embedding through the $H$-forcing, using the
master condition $j\image H$, and the fact that $j(\Q)$ is
${\lt}j(\k)$-directed closed.  Adding a further generic $j(H)\of
j(\Q)$, I lift to $j:V[G][H]\to M[j(G)][j(H)]$.  Again, using
$j\image\l$ as a seed, I generate a normal fine measure $\mu$ on
$P_\k\l$. The measure $\mu$ cannot have been added by $\Gtail$ or by
$j(H)$, and so $\mu$ is in $V[G][H]$, witnessing that $\k$ is
$\l$-supercompact there.\qed

\lemma. If $\phi(\Q,\k,A,G)$ fails in $V[G]$, then, below a condition,
$\Q$ destroys the supercompactness of $\k$.

\proof Let me prove the contrapositive. Suppose that every generic
extension by $\Q$ preserves the supercompactness of $\k$. Since $\phi$
respects the equivalence of forcing, I have merely to show that for
densely many $b\in\Q$ the relation $\phi(\Q_b,\k,A,G)$ holds in $V[G]$.
Fix any $b'\in\Q$; I intend to find a $b\lte b'$ such that
$\phi(\Q_b,\k,A,G)$ holds in $V[G]$.
Suppose that $H\of\Q$ is $V[G]$-generic below $b'$, where $\Q$ is
coded by some subset of $\l$. Suppose that $\k$ is still
$\theta$-supercompact in
$V[G][H]$, where $\theta\geq\l$ is such that $f$ represents $A$ with
respect to a $\theta$-supercompact embedding $j:V[G][H]\to M[j(G)][j(H)]$.
(By the Going-Down Lemma \GoingDownLemma, I may in fact take $\theta=\l$ here.)
It follows that $H\in M[j(G)][j(H)]$.  Factor $j(\P)$ as
$\P*\Qtilde*\Ptail$ and $j(G)$ as $G*g*\Gtail$, where  $g\of\Qtilde$ is
the (possibly trivial) stage $\k$ forcing in $j(\P)$.  Since $h$ jumps
over $\theta$, I know that the next forcing cannot occur until past $\theta$,
so $\Ptail$ is ${\lte}\theta$-closed.  Thus, $H\in M[G][g]$, and so, by the
Gap Forcing Theorem \GapForcingTheorem, since $M\of V$ it follows also that
$H\in V[G][g]$.  But $g\in V[G][H]$, and so
$V[G][g]=V[G][H]$. Thus, the forcing $\Qtilde$ and the forcing $\Q$
produce the same generic extension over $V[G]$. It follows that
$\RO(\Q_b)\iso\RO(\Qtilde_c)$ for some conditions $b\in\Q$ and
$c\in\Qtilde$. I may assume $b\lte b'$. Observe that
$\phi(\Qtilde,\k,A,G)$ holds in $M[G]$. By the
Gap Forcing Theorem, $M[G]$ and $V[G]$ have the
same $H(\theta\plus)$. Therefore, since $\phi$ is local,
$\phi(\Qtilde,\k,A,G)$ also holds in $V[G]$.
It follows, since $\phi$ respects the equivalence of forcing, that
$\phi(\Qtilde_c,\k,A,G)$ holds in $V[G]$, and hence also that
$\phi(\Q_b,\k,A,G)$ holds there, as desired.\qed

Thus, $V[G]$ is as required.\qed

\remark Remark on Closure. Again let me point out that directed closure is
only needed on the preservation side, to find a master condition. On the
destruction side, it is enough to assume that $\Q$ is $\ltk$-strategically
closed.

\section Epilogue: Fragility $\perp$ Superdestructibility

At first glance, fragility and superdestructibility seem to be made of
the same delicate material. But this is not so. In this epilogue, I
will show that neither property implies the other, and I will
construct models
which exhibit each of the four possibilities.

The notion of fragility first appeared in my first paper \cite[Ham94a],
and subsequently in my dissertation \cite[Ham94b], where I defined that
a large cardinal $\k$
is {\df fragile} when any forcing which preserves $\k\plus$ and $2^\ltk$
and adds a subset to $\k$ destroys the measurability of $\k$. The
notion of superdestructibility appeared first in \cite[Ham97b], where I
defined that a large cardinal $\k$ is {\df superdestructible} when any
$\ltk$-closed forcing which adds a subset of $\k$ destroys the
measurability of $\k$. I will now prove that these notions, though
similar, are actually independent.

\theorem Fragility $\perp$ Superdestructibility Theorem. Suppose $\k$
is a supercompact cardinal in $V$. Then in various forcing extensions where
$\k$ remains supercompact,
\points 1. $\k$ is both fragile and superdestructible.\cr
    2. $\k$ is fragile, but not superdestructible.\cr
        3. $\k$ is superdestructible, but not fragile.\cr
        4. $\k$ is neither fragile nor superdestructible.\cr

I will prove each of the four possibilities separately. So, for the remainder
of this paper, assume that $\k$ is a supercompact cardinal in $V$.  For
the first possibility, I will also show the surprising fact that one
can obtain a model in which $\k$ is fragile, superdestructible,
and, simultaneously, indestructible above $\k$---the supercompactness
of $\k$ is preserved by any $\ltk$-directed closed poset which adds no
subsets to $\k$.

\quiet\theorem Possibility One. There is a forcing extension in which
$\k$ is fragile, superdestructible, and, simultaneously,
indestructible above $\k$.

\proof In fact, the fragile measurability models of \cite[Ham94a] also
have superdestructibility. In order to get indestructibility above
$\k$, I will introduce here a wrinkle to the construction in the
Fragile Measurability Theorem 3.12 of \cite[Ham94a]. While familiarity
with that argument will ease comprehension of this one, I aim to give
here a complete, if terse, presentation.

So, suppose that $\k$ is supercompact in $V$. I may assume, by forcing
if necessary, that $V\satisfies\GCH$. Let $\ell$ be a Laver function
for $\k$. Our forcing $\P_{\k+1}=\P_\k*\Q_\k$ will be a reverse Easton
$(\k+1)$-iteration with nontrivial forcing only at inaccessible stages
$\g\in\dom(\ell)$. The forcing $\Q_\g$ at stage $\g$ will be one of three types.
First, the Laver function $\ell$ might instruct us to perform what I
will call the fragility forcing at stage $\g$. In this case, $\ell(\g)$
will hand us a pair $\<a_\g,d_\g>$, called {\df $\g$-data
packets} in \cite[Ham94a], such that $a_\g=\seq<\tlt_\g^\a\st\a<\g\plus>$
enumerates some relations on $\g$ such that $\ot\<\g,\tlt_\g^\a>=\a$ for
every $\a<\g\plus$,
and $d_\g=\seq<D_\g^\a\st\a<\g\plus>$ enumerates $P(\g)^V$. If there is
fragility forcing at some stage $\d<\g$, then it will have added a
sequence $\seq<C_\d^\a\st\a<\d\plus>$ of club subsets of $\d$, so I may
refer to these clubs when defining $\Q_\g$. Let $\Q_\g^\a$ be the poset
which adds, with conditions which are initial segments, a club set
$C_\g^\a\of\g$ with the property that if $\d$ is an inaccessible
cluster point of $C_\g^\a$, then $\<\g,\a>$ {\df reflects to}
$\<\d,\a'>$ for some $\a'$ in the sense that first of all there was
fragility forcing at stage $\d$, but secondly the data packets and
clubs agree: $$\tlt_\g^\a\restrict\d=\tlt_\d^{\a'}\qquad
D_\g^\a\intersect \d=D_\d^{\a'}\qquad C_\g^\a\intersect
\d=C_\d^{\a'}.$$ Let $\Q_\g$ be the $\ltg$-support product $\prod
\Q_\g^\a$. By a $\Delta$-system argument, this is $\g\plus$-c.c.  This
defines the fragility type forcing at stage $\g$.

The Laver function $\ell$, secondly, may instruct us to perform
Laver preparation forcing at stage $\g$. In this case, $\ell(\g)$ will
hand us a $\ltg$-directed closed poset $\Q_\g$, which will be our stage
$\g$ forcing provided that, additionally, it adds no new subsets to
$\g$.

Finally, third, the Laver function may instruct us to perform both
kinds of forcing. In this case, $\ell(\g)$ will hand us both the
$\g$-data packet, and also the (name of) a $\ltg$-directed closed poset
which adds no subsets to $\g$.  The stage $\g$ forcing will consist
of first performing the fragility forcing and then the
indestructibility forcing.

This completely describes the iteration $\P_\k$. The forcing $\Q_\k$ at
stage $\k$ will be of the fragility forcing type, with some specific
$\k$-data packet $\<a_\k,d_\k>$. Let $G*g\of\P_\k*\Q_\k$ be $V$-generic
for this forcing. Thus $g=\<C^\a_\k\st\a<\k\plus>$ is a sequence of club
subsets of $\k$ with the reflection property.
Let me now prove that $V[G][g]$ has the desired
properties. First observe that if there is forcing at stage $\g$, then
$\g\plus$ is not collapsed.

Let me now prove that $\k$ is fragile in $V[G][g]$.  Suppose towards a
contradiction that $H\of\Q$ is generic, adds a set $A\of\k$, preserves
$\k^{\ltk}$ and $\k\plus$, but that $\k$ remains measurable in
$V[G][g][H]$. Thus, there is an embedding $j:V[G][g][H]\to
M[j(G)][j(g)][j(H)]$. Since $\Q$ is not necessarily closed, the Gap
Forcing Theorem \GapForcingTheorem\ does not apply, but the set $A$
must be in $M[j(G)]$ since it can be added by neither $j(g)$ nor
$j(H)$.  Also, since $\k$ is necessarily an inaccessible closure point
of the clubs in $j(g)$, it follows that $j(G)$ must have performed
fragility forcing at stage $\k$. Moreover, by the reflection property
of the clubs in $j(g)$, since $\k$ is an inaccessible cluster point of
$j(C_\k^\a)$, it follows that the clubs added by the forcing at stage
$\k$ in $j(G)$ are $j(C_\k^\a)\intersect \k=C_\k^\a$. What is more,
since $j(\tlt_\k^\a)\restrict\k=\tlt_\k^\a$, it follows that $C_\k^\a$
is the generic used in the $\a^\th$ coordinate in the fragility forcing
at stage $\k$ in $j(\P)$. Since $\k\plus$ is preserved, this means that
the stage $\k$ fragility forcing in $j(G)$ is actually $g$.  Thus, it
follows that $j(G)=G*g*h*\gtail$, where $h$ is the stage $\k$
indestructibility forcing, if it exists, or $j(G)=G*g*\gtail$, if there
is no stage $\k$ indestructibility forcing.  In either case, since the
indestructibility forcing $h$ is not allowed to add subsets to $\k$,
the set $A$ must be in $M[G][g]$. Since $A$ is a subset of $\k$, I know
moreover that $A\in M[G][g\restrict\a]$ for some $\a<\k\plus$, and so
$A=\dot A_{G*g\restrict\a}$, for some name $\dot A\in M$.  Because of
the reflection $j(D^\a_\k)\intersect\k=(D_\k^\a)^M$ I know that
$P(\k)^M=P(\k)^V$, since these are both enumerated by $d_\k$, and thus
$\dot A\in V$.  Hence, $A\in V[G][g]$ contrary to our assumption that
$A$ was new.

Next, I will prove that $\k$ is superdestructible in $V[G][g]$.
Suppose that $H\of\Q$ is generic, where $\Q$ is $\ltk$-closed and $H$
adds a new set $A\of\k$, but that $\k$ is still measurable in
$V[G][g][H]$. Thus, there is an embedding $j:V[G][g][H]\to
M[j(G)][j(g)][j(H)]$. As above, I know that $A\in M[j(G)]$ by closure
considerations. I do not, however, know so easily that the stage $\k$
forcing of $j(G)$ is $g$, since it may be that $\Q$ collapsed
$\k\plus$. Nevertheless, since each $\Q_\g$ has dense sets as closed as
you like up to $\g$, it follows that $\P*\Q$ admits a gap below $\k$,
and thus, by the Gap Forcing Theorem \GapForcingTheorem,
$P(\k)^M=P(\k)^V$, and consequently ${\k\plus}^M={\k\plus}^V$. In fact,
$${\k\plus}^V={\k\plus}^M={\k\plus}^{M[j(G)]}={\k\plus}^{M[j(G)][j(g)][j(H)]}=
{\k\plus}^{V[G][g][H]}.$$ The first equality holds as I explained just
now. The second holds by my remark that the successor cardinals of
nontrivial forcing stages are not collapsed.
The third holds by the closure of the forcing. And the fourth
holds by the closure of the embedding. Thus, in fact, $\Q$ did not
collapse $\k\plus$. Since $\Q$ also preserves $2^\ltk$ and adds a subset to
$\k$, it follows from the already established
fact that $\k$ is fragile that $\Q$ destroys the measurability of $\k$.

Finally, I will argue that $\k$ is indestructible above $\k$ in $V[G][g]$.
This argument will also establish that $\k$ is supercompact in
$V[G][g]$. Suppose $H\of\Q$ is generic, where $\Q$ is $\ltk$-directed
closed and does not add new subsets to $\k$. Fix any $\l$ and let
$\theta\muchgt\l$.  Fix $j:V\to M$ a $\theta$-supercompact embedding in
$V$ such that $j(\ell)(\k)$ instructs us to first perform the composite
forcing at stage $\k$---first the fragility forcing and then $\Q$---and
such that $\dom(j(\ell))\intersect (\k,\theta]=\emptyset$. I will lift $j$ to
$V[G][g][H]$. First, I can lift to $j:V[G]\to M[j(G)]$ by using $g*H$
as the stage $\k$ generic, and then forcing to add a tail $\Gtail$.
Now I have
to construct a master condition below $j\image g$ as in \cite[Ham94a],
and force below it to lift to $j:V[G][g]\to M[j(G)][j(g)]$.  The master
condition is simply the condition $p\in j(\Q_\k)$, with support
$j\image\k\plus$, such that $p(j(\a))=\bar C_\k^\a$, where $\bar
C_\k^\a=C_\k^\a\union\{\k\}$. Below this condition there is in
$j(\Q_\k)$ a dense set which is ${\lte}\theta$-closed, namely, the set of
conditions $q$ which mention a point above $\theta$ on every coordinate
in their support. So I can, by forcing over $j(\Q_\k)$ below the master
condition, lift $j$ through the $\Q_\k$ forcing.  Next, use the
directed closure of $j(\Q)$ to find a master condition below $j\image
H$, and lift fully to $j:V[G][g][H]\to M[j(G)][j(g)][j(H)]$. This
embedding lives in $V[G][g][H][\Gtail][j(g)][j(H)]$. But using
$j\image\l$ as a seed, I conclude that there is a measure witnessing
$\l$-supercompactness which could not have been added by the tail
forcing $\Gtail*j(g)*j(H)$. So the measure lives in $V[G][g][H]$, and so
$\k$ is $\l$-supercompact there, as desired.
In the case that $\Q$ is trivial, we conclude
also that $\k$ is supercompact in $V[G][g]$.
This completes the proof.\qed Possibility One

\quiet\theorem Possibility Two. There is a forcing extension in which
$\k$ is fragile, but not superdestructible. In fact, $\k$ can be made
simultaneously fragile and indestructible by any $\ltk$-directed closed
forcing which collapses $\k\plus$.

\proof Here I will modify the previous argument. Again I will perform a
$\k+1$ reverse Easton iteration $\P_\k*\Q_\k$, where at each stage $\g$
I perform one of three kinds of forcing. First, the Laver function may
instruct us, as before, to perform the fragility forcing. Second, the
Laver function may instruct us to perform $\ltg$-directed closed
forcing $\Q_\g$, and we will oblige, {\it provided} that $\Q_\g$
collapses $\g\plus$. Finally, third, the Laver function may instruct us
to perform both of the previous two types of forcing. Let $G*g$ be
$V$-generic for $\P_\k*\Q_\k$, where $\Q_\k$ is as previously the
stage $\k$ fragility forcing using the $\k$-data packet $\<a_\k,d_\k>$,
and let me show that $V[G][g]$ has the
properties that we seek.

First, I will prove that $\k$ is supercompact in $V[G][g]$.
Fix any $\l$, and select $\theta\muchgt\l$, and a $\theta$-supercompact
embedding $j:V\to M$ such that $j(\ell)(\k)$ tells us to perform just the
fragility forcing, using $\<a_\k,d_\k>$, and that
$\dom(j(\ell))\intersect (\k,\theta]=\emptyset$.
Thus, $j(\P)=\P*\Q_\k*\Ptail$, where $\Ptail$ is
${\lte}\theta$-closed. As usual, force over the tail, and lift the embedding
to $j:V[G]\to M[j(G)]$, where $j(G)=G*g*\Gtail$. Now use the master condition
argument from Possibility One to lift through the $\Q_\k$ forcing. This gives
$j:V[G][g]\to M[j(G)][j(g)]$ in $V[G][g][\Gtail][j(g)]$. The measure on
$P_\k\l$ germinated by the seed $j\image\l$ must lie in $V[G][g]$, so $\k$ is
$\l$-supercompact there.

Next, I will establish that $\k$ is fragile in $V[G][g]$.
This is nearly identical to the corresponding argument in Possibility
One. Suppose that $\k$ is measurable in $V[G][g][H]$, where $H\of\Q$ adds a
new subset $A\of\k$, but preserves $\k^{\ltk}$ and $\k\plus$.  Thus,
there is an embedding $j:V[G][g][H]\to M[j(G)][j(g)][j(H)]$. Moreover,
since $\k$ is an inaccessible cluster point of the club sets in $j(g)$, namely
$j(C_\k^\a)$, it follows that the $\a^\th$ generic club added at stage $\k$ is
$j(C_\k^\a)\intersect \k=C_\k^\a$. Thus, since as before $\k\plus$ is preserved,
$j(G)=G*g*\Gtail$ for some $\Gtail$.
Furthermore, by the reflection property of
the stage $j(\k)$ data packets to the stage $\k$ data packets, I also know
$P(\k)^M=P(\k)^V$, since both are enumerated by $d_\k$.
Also, $\k\plus$ is the same in all the models since
it is not collapsed by $G$, $g$, or $H$. Thus there was no supplementary
Laver forcing at stage $\k$ in $j(\P)$. So I can proceed as before,
and obtain the contradiction involving $\dot A$.

Finally, let me prove that the supercompactness of $\k$ is preserved over
$V[G][g]$ by any $\ltk$-directed closed poset which collapses
$\k\plus$.
Suppose that $H\of\Q$ is generic, where $\Q$ is $\ltk$-directed
closed and collapses $\k\plus$. Fix any $\l$ and pick
$\theta\muchgt\l$. Fix $j:V\to M$ a $\theta$-supercompact embedding
such that $j(\ell)(\k)$ instructs us to perform the fragility forcing with
$\<a_\k,d_\k>$,
followed by $\Q$.  Also I will need that $\dom(j(\ell))\intersect
(\k,\theta]=\emptyset$.  Let me proceed to lift $j$ to the forcing
extension. The stage $\k$ forcing in $j(\P_\k)$ is no problem since I have
$g*H$. Force over the tail to get $\Gtail$, and lift to $j:V[G]\to M[j(G)]$,
where $j(G)=G*g*H*\Gtail$.  Now use the master condition below $j\image
g$ and force to add $j(g)$, lifting to $j:V[G][g]\to M[j(G)][j(g)]$.
Similarly, using the directed closure, I can find a master condition
below $j\image H$ and lift to $j:V[G][g][H]\to M[j(G)][j(g)][j(H)]$.
Finally, use $j\image \l$ as a seed and observe that $\k$ is still
$\l$-supercompact in $V[G][g][H]$.  The measure could not have been
added by $\Gtail*j(g)*j(H)$, because of closure, and therefore lies in
$V[G][g][H]$. This completes the proof.\qed Possibility Two

\quiet\theorem Possibility Three. There is a forcing extension in which
$\k$ is superdestructible, but not fragile.

\proof First force with the Laver preparation to $V[G]$ where $\k$ is
indestructible. Then, perform any small forcing $g\of\Q$, such as adding a
single Cohen real. It follows by
the Superdestruction Theorem of \cite[Ham97b] that $\k$ is
superdestructible in $V[G][g]$, and I will now show that $\k$ is not
fragile there. For this, it suffices to show that forcing over
$\R=\add(\k,1)^{V[G]}$ will preserve the supercompactness of $\k$ over
$V[G][g]$. Certainly forcing with $\R$ over $V[G]$ preserves the
supercompactness of $\k$, since $\k$ is indestructible in $V[G]$. Also,
the further small forcing
$\Q$ preserves supercompactness again. So over $V[G]$, the forcing
$\R\cross \Q$ preserves supercompactness. Thus, by rearranging the order of
the forcing, $\Q\cross \R$ also preserves supercompactness. Thus,
$\R$ preserves supercompactness over $V[G][g]$, as desired.\qed Possibility
Three

\quiet\theorem Possibility Four. There is a forcing extension in which
$\k$ is neither fragile nor superdestructible.

\proof The Laver preparation makes $\k$ indestructible, and therefore
neither fragile nor superdestructible.\qed Possibility Four

\section Bibliography

\tenpoint

\nopagenumbers
\parindent=0pt
\newbox\Article
\newbox\Journal
\newbox\Author
\newbox\Vol
\newbox\No
\newbox\Year
\newbox\Page
\newbox\Book
\newbox\Publisher
\newbox\Pubaddr
\newbox\Key
\newbox\Editor
\newbox\Comment
\newbox\Note
\def\entry#1#2\par{\item{#1\quad}\hskip-1.1em#2\par}
\def\article#1{\setbox\Article=\hbox{\sl #1, }}
\def\journal#1{\setbox\Journal=\hbox{\rm #1 }}
\def\author#1{\setbox\Author=\hbox{\sc #1, }}
\def\vol#1{\setbox\Vol=\hbox{\bf #1 }}
\def\no#1{\setbox\No=\hbox{no. #1 }}
\def\year#1{\setbox\Year=\hbox{\rm({\oldstyle #1}) }}
\def\page#1{\setbox\Page=\hbox{\rm p. #1 }}
\def\book#1{\setbox\Book=\hbox{\it #1, }}
\def\publisher#1{\setbox\Publisher=\hbox{\rm #1, }}
\def\pubaddr#1{\setbox\Pubaddr=\hbox{\rm #1, }}
\def\key#1{\setbox\Key=\hbox{#1}}
\def\editor#1{\setbox\Editor=\hbox{\rm(#1, Ed.) }}
\def\comment#1{\setbox\Comment=\hbox{\rm #1}}
\def\note#1{\setbox\Note=\hbox{\rm #1 }}
\def\ref#1\par{\smallskip{#1
  \entry{\ifhbox\Key\unhbox\Key\else[\ ]\fi}%
  \unhbox\Author\unhbox\Note
  \ifhbox\Book \unhbox\Book\unhbox\Publisher\unhbox\Pubaddr
               \unhbox\Editor\unhbox\Page\unhbox\Year\unhbox\Comment
  \else \unhbox\Article\unhbox\Journal\unhbox\Vol\unhbox\No\unhbox\Editor
        \unhbox\Page\unhbox\Year\unhbox\Comment\fi\par}}

\ref
\author{Arthur W. Apter}
\article{Laver Indestructibility and the Class of Compact Cardinals}
\journal{Journal of Symbolic Logic}
\comment{(to appear)}
\key{[Apt96]}

\ref
\author{Joel David Hamkins}
\article{Fragile Measurability}
\journal{Journal of Symbolic Logic}
\year{1994}
\vol{59}
\page{262-282}
\key{[Ham94a]}

\ref
\author{Joel David Hamkins}
\article{Fragile Measurability; Lifting and Extending Measures}
\year{1994}
\comment{UC Berkeley dissertation}
\key{[Ham94b]}

\ref
\author{Joel David Hamkins}
\article{Canonical Seeds and Prikry Trees}
\journal{Journal of Symbolic Logic}
\vol{62}
\no{2}
\year{1997}
\key{[Ham97a]}

\ref
\author{Joel David Hamkins}
\article{Small Forcing Makes Any Cardinal Superdestructible}
\journal{Journal of Symbolic Logic}
\vol{62}
\no{4}
\year{1997}
\comment{(to appear)}
\key{[Ham97b]}

\ref
\article{Superdestructibility: A Dual to Laver
Indestructibility}
\author{Joel David Hamkins \& Saharon Shelah}
\journal{Journal of Symbolic Logic}
\comment{(to appear)}
\key{[HamShl]}

\ref
\author{Akihiro Kanamori}
\book{The Higher Infinite}
\publisher{Springer Verlag}
\year{1994}
\key{[Kan94]}

\ref
\author{Kimchi \& Magidor}
\article{The Independence between the Concepts of Compactness and
Supercompactness}
\comment{circulated manuscript}
\key{[KimMag]}

\ref
\author{Richard Laver}
\article{ Making the Supercompactness of $\kappa$ Indestructible Under
 $\kappa$-Directed Closed Forcing}
\journal{Israel Journal Math}
\vol{29}
\year{1978}
\page{385-388}
\key{[Lav78]}

\ref
\author{Silver, Jack H.}
\article{On the Singular Cardinals Problem}
\journal{Proceedings International Congress of Mathematics Vancouver}
\year{1974}
\page{265-268}
\key{[Sil74]}

\ref
\author{Solovay, Robert M.}
\article{Strongly Compact Cardinals and the \GCH}
\journal{Proceedings of the Tarski Symposium, Proceedings of
Symposia in Pure Mathematics}
\year{1974}
\vol{25}
\page{365-372}
\key{[Sol74]}

\ref
\author{W. Hugh Woodin}
\article{A Supercompact Cardinal Whose Weak Compactness is Destroyed by
         $\add(\k,1)$}
\comment{(personal communication)}
\key{[W]}

\bye